\documentclass[10pt]{article}
\usepackage{amsfonts}
\usepackage{amsfonts,latexsym,amsmath,mathrsfs,amscd,geometry}
\usepackage{amssymb}
\usepackage{latexsym}
\usepackage{color}
\usepackage{upgreek}
\usepackage{amsmath}
\numberwithin{equation}{section}

\newcommand \nc{\newcommand}
\newtheorem{theorem}{Theorem}[section]
\newtheorem{lemma}[theorem]{Lemma}
\newtheorem{proposition}[theorem]{Proposition}

\newtheorem{remark}[theorem]{Remark}

\nc{\ba}{\begin{array}}\nc{\ea}{\end{array}}
\nc{\be}{\begin{eqnarray}}\nc{\ee}{\end{eqnarray}}
\nc{\beq}{\begin{equation}}\nc{\eeq}{\end{equation}}
\nc{\bex}{\begin{eqnarray*}}\nc{\eex}{\end{eqnarray*}}
\nc{\btm}{\begin{theorem}} \nc{\etm}{\end{theorem}}
\nc{\blm}{\begin{lemma}} \nc{\elm}{\end{lemma}}
\nc{\R}{\mathbb{R}} \nc{\va}{\nu} \nc{\ls}{\limits}

\def\pf{\noindent{\bf Proof.\quad}}\def\endpf{\hfill$\Box$}

\allowdisplaybreaks

\usepackage[colorlinks,linkcolor=blue]{hyperref}

\begin{document}
\title{Incompressible limit of a compressible Oldroyd-B model}
\author{Sili Liu\footnote{School of Mathematics and Statistics, Changsha University of Science and Technology, Changsha 410114, China. E-mail: slliu@csust.edu.cn} $^\dag$,\, Yingshan Chen\footnote{(Corresponding Author) School of Mathematics, South China University of Technology, Guangzhou 510631, China. E-mail: mayshchen@scut.edu.cn.}}

\maketitle

\begin{abstract}
In this paper, we consider the Cauchy problem for a compressible Oldroyd-B model in three dimensions. Under some smallness assumptions on the initial data, we obtain the global wellposedness of strong solution with uniform regularity. Moreover, we provide a rigorous justification for the link between the compressible model and the corresponding incompressible model via incompressible limits.

\end{abstract}

\bigbreak \textbf{{\bf Key Words}:}  Compressible Oldroyd-B model; Incompressible limit; Convergence rate.

\bigbreak  {\textbf{AMS Subject Classification 2020:} 76N10, 76N17, 74H40.

\section{Introduction}

In this paper, we study a compressible Oldroyd-B model considered in a series of papers such as \cite{BLS, LLW, LWW, LZ}. Since our main concern is the incompressible limit, we write the model into the following dimensionless form in $\mathbb{R}^3\times\mathbb{R}^+$:
\beq\label{clc-01}
\left\{
\begin{split}
&\rho_t+\mathrm{div} (\rho u)=0,\\
&(\rho u)_t+\mathrm{div} (\rho u\otimes u)+\frac{1}{\epsilon^2}\nabla P(\rho)=\mu_1\Delta u+\mu_2\nabla \mathrm{div}u
+\mathrm{div} \bigg(\frac{\beta}{k}\mathbb{T}-(\beta L\eta+\bar{\mathfrak{z}}\eta^{2})\mathbb{I}\bigg),\\
&\eta_t+\mathrm{div} (\eta u)=\nu\Delta \eta,\\
&\mathbb{T}_t+\mathrm{div} (u\mathbb{T})-(\nabla u
\mathbb{T}+\mathbb{T}\nabla^Tu)
=\nu\Delta \mathbb{T}+\frac{k A_0}{2}\eta\mathbb{I}-\frac{A_0}{2}\mathbb{T},
\end{split}
\right.
\eeq
where the pressure $P(\rho)$ and the density $\rho(x,t)\geq0$ of the fluid are supposed to be related by the typical power law relation for simplicity:
\[P(\rho)=a\rho^{\gamma}\]
for some known constants $a>0$, $\gamma>1$; $u(x,t)\in\mathbb{R}^{3}$ denotes the velocity field of the fluid. $\mu_1>0$ and $\mu_2>0$ are viscosity coefficients. The polymer number density $\eta(x,t)\geq0$ represents the integral of
the probability density function $\psi$ which is a microscopic variable in the modelling of dilute polymer chains, i.e.,
\[\eta(x,t)=\int_{\mathbb{R}^3}\psi(x,q,t)\,\mathrm{d}q,\]
where $\psi$ is governed by the Fokker-Planck equation. The extra stress tensor $\mathbb{T}(x,t)= (\mathbb{T}_{i,j})(x,t)\in\mathbb{R}^{3\times 3},\ 1\leq i,j\leq 3$ is a positive definite symmetric matrix, and the notation $\mathrm{div} (u\mathbb{T})$ is understood as
\begin{equation*}
\left(\mathrm{div} (u\mathbb{T})\right)_{i,j} = \mathrm{div} (u\mathbb{T}_{i,j}), \quad 1\leq i,j \leq 3.
\end{equation*}
The constant parameter $\epsilon>0$ is the Mach number, $\beta=\frac{n_p}{n_p+n_s}$, with $n_p$ signifying the polymeric viscosity and $n_s$ the viscosity of the solvent, $\nu$ is the centre-of-mass diffusion coefficient and other parameters $k, L, \bar{\mathfrak{z}}, A_0$ are all positive and known numbers, whose meanings will be explained in detail when we (formally) derive system \eqref{clc-01} in the following subsection.

It is known that the compressible Oldroyd-B model can be derived as a macroscopic closure of Navier-Stokes-Fokker-Planck system, which is a micro-macro model
describing dilute polymeric fluids, in Hookean bead-spring chain setting (see \cite{BLS}). In the following we shall use the idea in \cite{BLS} to (formally) derive \eqref{clc-01} from the dimensionless form of the Navier-Stokes-Fokker-Planck system considered in \cite{SW}.

\subsection{Formal derivation of the system $(\ref{clc-01})$ }
In order to derive system $(\ref{clc-01})$, we recall the dimensionless form of the Navier-Stokes-Fokker-Planck system considered in \cite{SW}:
\beq\label{f-01}
\left\{
\begin{split}
&\rho_t+\mathrm{div}_x (\rho u)=0,\\
&(\rho u)_t+\mathrm{div}_x (\rho u\otimes u)+\frac{1}{\text{Ma}^2}\nabla_x P(\rho)=\frac{1}{\text{Re}}\mathrm{div}_x \mathbf{S}(u)+\frac{1}{\text{Fr}^2}\rho f
+\frac{1}{\text{Re}}\mathrm{div}_x \tau_1 -\frac{\tilde{\mathfrak{z}}}{\text{Ma}^2}\nabla_x(\eta^{2}),\\
&\psi_t+\mathrm{div}_x (u\psi)+\sum_{i=1}^K \mathrm{div}_{q_i}\big((\nabla_x u)q_i\psi\big)
=\nu\Delta_x \psi+\frac{1 }{4\text{De}}\sum_{i=1}^K \sum_{j=1}^K A_{i,j}\mathrm{div}_{q_i}\bigg(M\nabla_{q_i}\big(\frac{\psi}{M}\big)\bigg),
\end{split}
\right.
\eeq
where $\text{Ma}$ is the Mach number, $\text{Re}$ is the Reynolds number, $\text{Fr}$ is the Froude number, and $\text{De}$ is the Deborah number. The constant matrix $A=(A_{ij}), 1\leq i,j\leq K$, called the Rouse matrix, is symmetric and positive definite. We denote by $A_0$ the smallest eigenvalue of $A$; clearly, $A_0 > 0$. In addition, the Newtonian shear stress tensor $\mathbf{S}(u)$ is defined by
 \[\mathbf{S}(u)=\mu^S\bigg(\nabla_xu+\nabla_x^Tu
 -\frac{2}{3}(\mathrm{div}_xu)\mathbb{I}\bigg)+\mu^B(\mathrm{div}_xu)\mathbb{I}.\]
A direct calculation gives
\begin{align}\label{f-09}
\mathrm{div}_x\mathbf{S}(u)=\mu^S\Delta_xu
+\Big(\mu^B+\frac{\mu^S}{3}\Big)\nabla_x\mathrm{div}_xu
=\mu_1\Delta_xu+\mu_2\nabla_x\mathrm{div}_xu,
\end{align}
with $\mu_1=\mu^S>0$ and $\mu_2=\mu^B+\frac{\mu^S}{3}>0$.

In a bead-spring chain model consisting of $L(=K+1)$ beads coupled with $K$ elastic springs representing a polymer chain, the probability density function $\psi$ depends not only on $x$ and $t$, but also on the conformation vector $q=(q_1^T,\cdots,q_K^T)^T\in D:=D_1\times\cdots D_K\subset\mathbb{R}^{3K}$, with $q_i$ representing the 3-component conformation/orientation vector of the $i$th spring in the chain. Typically $D_i$ is the whole space $\mathbb{R}^{3}$ or a bounded open ball centered at the origin $0 $ in $\mathbb{R}^{3}$, for each $i=1,\ldots,K$. Where $K=1$, the model is referred to as the dumbbell model. Here we consider the Hookean bead--spring chain model, where $D_i=\mathbb{R}^{3}$ for all $i\in\{1,\ldots,K\}$, and the elastic spring--force $F_i:q_i \in D_i\mapsto U'_i(\frac{1}{2}|q_i|^2)q_i\in\mathbb{R}^{3}$ and the spring potential $U_i:\mathbb{R}_{\geq0}\mapsto\mathbb{R}_{\geq0}$ of the $i$th spring in the chain are defined by
\begin{align}\label{f-02}
F_i(q_i)=q_i\;\;\text{for all}\;\;q_i\in D_i,\;\; U_i(s)=s\;\;\text{for all}\;\; s\geq0,\;\;i=1,\ldots K.
\end{align}
On the right-hand side of \eqref{f-01}$_2$, the 3-component vector function $f$ is the nondimensional density of body forces and the elastic extra-stress tensor is of the form:
\begin{align*}
\tilde{\tau}(\psi)(x,t)&:=\frac{1}{\text{Re}}\tau_1(\psi)(x,t)
-\frac{1}{\text{Ma}^2}\Big(\int_{D\times D}\tilde{\gamma}(q,q')\psi(x,q,t)\psi(x,q',t)\;\mathrm{d}q\mathrm{d}q'\Big)\mathbb{I}
\\
&=\frac{1}{\text{Re}}\tau_1(\psi)(x,t)
-\frac{\tilde{\mathfrak{z}}}{\text{Ma}^2}\eta^2(x,t)\mathbb{I},
\end{align*}
where $\tilde{\gamma}:D\times D\rightarrow \mathbb{R}_{\geq0}$ is a smooth, $t$--independent, $x$--independent and $\psi$--independent interaction kernel, which we take here for the sake of simplicity to be
\[\tilde{\gamma}(q,q')=\tilde{\mathfrak{z}},
\;\;\text{where}\;\;\tilde{\mathfrak{z}}\in\mathbb{R}_{\geq0}.\]
Moreover, $\tau_1$ is the Kramers expression:
\begin{align}\label{f-03}
\tau_1:= \frac{\beta}{\text{De}}\bigg(\sum_{i=1}^K\mathbf{C}_i(\psi)
-L\int_D\psi\;\mathrm{d}q\,\mathbb{I}\bigg),
\end{align}
where $\beta$ is a constant dependent on the polymeric viscosity and the viscosity of the solvent, and
\begin{align}\label{f-04}
\mathbf{C}_i(\psi)(x,t):
=\int_D\psi(x,q,t)U'_i(\frac{1}{2}|q_i|^2)q_iq_i^T\;\mathrm{d}q,\;\;
i=1,\ldots,K.
\end{align}
Hence, from \eqref{f-02} and \eqref{f-04}, \eqref{f-03} can be rewritten as in the Hookean case
\begin{align}\label{f-05}
\tau_1:= \frac{\beta}{\text{De}}\bigg(\frac{1}{k}\mathbb{T}(x,t)
-L\eta(x,t)\mathbb{I}\bigg),
\end{align}
with \[\mathbb{T}(x,t)=k\sum_{i=1}^K\int_D\psi(x,q,t)q_iq_i^T\;\mathrm{d}q,\]
where $k>0$ is the product of the Boltzmann constant and the absolute temperature.

Based on \cite{SW}, the expression appearing in the second term on the right-hand side of \eqref{f-01}$_3$ can be rewritten as follows:
\[M\nabla_{q_i}\big(\frac{\psi}{M}\big)=\nabla_{q_i}\psi+\psi q_i,\]
and
\begin{align}\label{f-06}
\psi|q_i|\rightarrow0,\;\;\nabla_{q_i}\psi\cdot\frac{q_i}{|q_i|}\rightarrow0,
\;\;\text{as}\;\;|q_i|\rightarrow\infty,\;\;\text{for all}\;\;
(x,t)\in\mathbb{R}^3\times(0,T],\;\;i=1,\ldots,K.
\end{align}

Finally, recalling the definition of $\eta(x,t)$, integrating the Fokker--Planck equation \eqref{f-01}$_3$ with respect to $q$ over $D$ and using \eqref{f-06}, one can follow \cite{BLS} step by step to deduce the equation of $\eta$:
\begin{align}\label{f-07}
\eta_t+\mathrm{div}_x(u\eta)=\nu\Delta_x\eta.
\end{align}
Then, considering the special case of the model $K=1$, multiplying the Fokker--Planck equation \eqref{f-01}$_3$ by $k\sum_{i=1}^Kq_iq_i^T$ and then integrating the result with respect to $q$ over $D$, it easily follows from the same technical arguments in Section 2.2 of \cite{BLS} that
\begin{align}\label{f-08}
\mathbb{T}_t+\mathrm{div}_x(u\mathbb{T})
-(\nabla_xu\mathbb{T}+\mathbb{T}\nabla_x^Tu)=\nu\Delta_x\mathbb{T}
+\frac{kA_0}{2\text{De}}\eta\mathbb{I}-\frac{A_0}{2\text{De}}\mathbb{T}.
\end{align}
Thus, \eqref{f-01} together with \eqref{f-09}, \eqref{f-05}, \eqref{f-07} and \eqref{f-08}, one has
\beq\label{f-10}
\left\{
\begin{split}
&\rho_t+\mathrm{div}_x (\rho u)=0,\\
&(\rho u)_t+\mathrm{div}_x (\rho u\otimes u)+\frac{1}{\text{Ma}^2}\nabla_x P(\rho)-\frac{\mu_1}{\text{Re}}\Delta_xu-\frac{\mu_2}{\text{Re}}\nabla_x\mathrm{div}_xu\\
&\qquad\qquad\quad=\frac{1}{\text{Fr}^2}\rho f
+\frac{1}{\text{Re}} \frac{\beta}{\text{De}}\mathrm{div}_x \bigg(\frac{1}{k}\mathbb{T}
-L\eta\mathbb{I}\bigg) -\frac{\tilde{\mathfrak{z}}}{\text{Ma}^2}\nabla_x(\eta^{2}),\\
&\eta_t+\mathrm{div}_x(u\eta)=\nu\Delta_x\eta,\\
&\mathbb{T}_t+\mathrm{div}_x(u\mathbb{T})
-(\nabla_xu\mathbb{T}+\mathbb{T}\nabla_x^Tu)=\nu\Delta_x\mathbb{T}
+\frac{kA_0}{2\text{De}}\eta\mathbb{I}-\frac{A_0}{2\text{De}}\mathbb{T}.
\end{split}
\right.
\eeq
Choosing $\text{Ma}=\epsilon$, $\text{Re}=\text{De}=1$, $\tilde{\mathfrak{z}}=\bar{\mathfrak{z}}\text{Ma}^2=\bar{\mathfrak{z}}\epsilon^2$ and $f=0$ as in \cite{SW}, we can get \eqref{clc-01} from \eqref{f-10} directly.

\subsection{Historical work and our contributions}
In this subsection, we recall some mathematical results for the Oldroyd-B models briefly. Although the mathematical study of Oldroyd-B models attracts a lot of attention, there are few results about the system \eqref{clc-01}. In \cite{DMPS}, we notice that the compressible Oldroyd-B model \eqref{clc-01} and some related non-isothermal models can be derived in the thermodynamic framework. This idea can be traced back to \cite{RS} and has been widely used in \cite{BMP, HMPST, MPSS, MRT, SPMR}. \cite{BFM} studied the relevant models and proved the global existence of weak solutions under the condition that the pressure term index $\gamma$ is large enough. When ignoring the term $\eta$ and assuming that the stress tensor diffusion coefficient is zero, Fang and Zi \cite{FZ} proved the local existence and uniqueness of the strong solution to a compressible Oldroyd-B model and established a blow-up criterion for strong solutions. The research on the incompressible limit of this simplified model can be seen in the seminal work by Lei (\cite{Lei}) where the uniform local regularity was derived, see \cite{FZ1} for the framework in Besov space. The compressible Oldroyd-B model \eqref{clc-01} was derived as a macroscopic closure of of the Navier-Stokes-Fokker-Planck system in \cite{BLS}, where the existence of global finite-energy weak solutions with arbitrarily large initial data in two dimensions was also proved. However, the existence of global solutions strong or weak with arbitrarily large initial data remains open in three dimensions, except the recent result \cite{Lu-Pokorny20} for a simplified model where the extra stress tensor is assumed to be a scalar function multiplied by a unit matrix. The local existence and uniqueness of strong solutions in two and three dimensions were proved by Lu and Zhang \cite{LZ}, where a blow-up criterion and the weak-strong uniqueness were also studied. For the case of small initial data, Wang and Wen \cite{W1W} showed the global well-posedness and optimal time decay of strong solutions in $H^3$ spaces. Recently, in \cite{LLW,LWW}, the first author and her collaborators investigated the vanishing of centre-of-mass diffusion and the inviscid case either for global well-posedness or for long time behavior. To the best of our knowledge, there is no result on the relationship between the compressible Oldroyd-B model \eqref{clc-01} and its incompressible counterpart.

Inspired by a recent work for the incompressible limits of Ericksen-Leslie hyperbolic liquid crystal model (\cite{GJLLT}), we consider the density $\rho$ with a small perturbation around the equilibrium state $1$ as $\rho^{\epsilon}=1+\epsilon\phi^{\epsilon}$, and denote $\tau^\epsilon=\mathbb{T^\epsilon}-k\eta^\epsilon \mathbb{I}$ as in \cite{LLW,LWW,Lu-Pokorny20}, then the system \eqref{clc-01} becomes
\beq\label{ma1}
\left\{
\begin{split}
&\partial_t\phi^{\epsilon}+u^{\epsilon}\cdot\nabla\phi^{\epsilon}
+\phi^{\epsilon}\mathrm{div}u^{\epsilon}+\frac{1}{\epsilon}\mathrm{div}u^{\epsilon}=0,\\
&\partial_t u^{\epsilon}+ u^{\epsilon}\cdot \nabla u^{\epsilon}+\frac{1}{\epsilon}\frac{P'(\rho^{\epsilon})}{\rho^{\epsilon}} \nabla\phi^{\epsilon} +\frac{1}{\rho^{\epsilon}}\nabla\big(\beta(L-1)\eta^{\epsilon}+\bar{\mathfrak{z}}(\eta^{\epsilon})^{2}\big)
=\frac{\mu_1}{\rho^{\epsilon}}\Delta u^{\epsilon}+\frac{\mu_2}{\rho^{\epsilon}}\nabla \mathrm{div}u^{\epsilon}+\frac{\beta}{k}\frac{1}{\rho^{\epsilon}}\mathrm{div} \tau^{\epsilon},\\
&\partial_t\eta^{\epsilon}+\mathrm{div} (\eta ^{\epsilon} u^{\epsilon})=\nu\Delta \eta^{\epsilon},\\
&\partial_t\tau^{\epsilon}+\mathrm{div} (u^{\epsilon}\tau^{\epsilon})-(\nabla u^{\epsilon}
\tau^{\epsilon}+\tau^{\epsilon}\nabla^Tu^{\epsilon})-k\eta^{\epsilon}(\nabla u^{\epsilon}+\nabla^Tu^{\epsilon})
=\nu\Delta \tau^{\epsilon}-\frac{A_0}{2}\tau^{\epsilon}.
\end{split}
\right.
\eeq

From mathematical point of view, it is reasonable to expect that, as $\rho^{\epsilon}\rightarrow1$, the first equation in \eqref{ma1}  yields the incompressible condition $\mathrm{div}u=0$. Moreover, suppose that $(u^{\epsilon}, \eta^{\epsilon},\tau^{\epsilon})\rightarrow(u,\eta,\tau)$ exists as  $\epsilon\rightarrow 0$, then (at least formally)  we obtain the following incompressible Oldroyd-B model:
\beq\label{limit1}
\left\{
\begin{split}
&\mathrm{div}u=0,\\
&\partial_t u+ u\cdot \nabla u-\mu_1\Delta u+\nabla\pi
=\frac{\beta}{k}\mathrm{div} \tau,\\
&\partial_t\eta+\mathrm{div} (\eta  u)=\nu\Delta \eta,\\
&\partial_t\tau+\mathrm{div} (u\tau)-(\nabla u
\tau+\tau\nabla^Tu)-k\eta(\nabla u+\nabla^Tu)
=\nu\Delta \tau-\frac{A_0}{2}\tau,
\end{split}
\right.
\eeq
where $\nabla\pi$ is the ``limit" of $\frac{1}{\epsilon^2} \nabla P(\rho^{\epsilon}) +\nabla\big(\beta(L-1)\eta^{\epsilon}+\bar{\mathfrak{z}}(\eta^{\epsilon})^{2}\big)$.

The main objective of this paper is to rigorously prove that the strong solution of \eqref{ma1} would converge to the solution of of \eqref{limit1}.

When $\epsilon$ is any a fixed positive constant, Wang and Wen \cite{W1W} conquered the technical difficulties due to the terms related to $\eta$ in the equation for the extra stress tensor and in the momentum equation, and proved the global existence and uniqueness of solutions to the system \eqref{clc-01} by deducing some key estimates for the polymer number density and its spatial derivatives. Compared to \cite{W1W}, our main difficulties are to deal with the singular terms $\frac{1}{\epsilon}\mathrm{div}u^\epsilon$ and $\frac{1}{\epsilon}\frac{P'(\rho^\epsilon)}{\rho^\epsilon}\nabla\phi^\epsilon$ in equation \eqref{ma1} and to derive some uniform $($in $\epsilon$$)$ estimates. 

The main results of this paper can be viewed as an extension of \cite{Lei} in the sense that the uniform regularity and incompressible limits are derived globally in time. On one hand, from the perspective of the equation structure, the coefficient of the fourth term on the left-hand side of \eqref{ma1}$_4$ is not a constant but rather satisfy a transport-diffusion equation, and two gradient terms are added to the momentum equation. On the other hand, we rigorously show that the incompressible limit of the global strong solution as well as the convergence rate of the limit. 
For more results on this topic in different models (Navier-Stokes, MHD, etc.), please refer to \cite{DHWZ,JJL,JJLX,LM1,T} and the references therein.}

\subsection{Notations and main results}
In order to simplify the presentation, we introduce some notations. For all $p\in[1,+\infty]$, $L^p:=L^p(\mathbb{R}^3)$ is the standard $L^p$ space. For
$p\in[1,+\infty)$ and some positive weight function $\varphi:\mathbb{R}^3\rightarrow \mathbb{R}^+$, the weighted space $L^p_\varphi:=L^p_\varphi(\mathbb{R}^3)$  endow with the norm $\|\cdot\|_{L^p_\varphi}=\|\cdot \varphi^\frac{1}{p}\|_{L^p}$. For integer $m\geq0$, $H_\varphi^m$ is the weighted Sobolev space with the norm
\[\|f\|_{H_\varphi^m}=\big(\sum\limits_{|\alpha|\leq m}\|\partial^\alpha f\|_{L^2_\varphi}^2\big)^\frac{1}{2}.\]
$\|f\|_{H_{\text{loc}}^m}<+\infty$ means that $\|f\|_{H^m(\Omega)}<+\infty$, for any compact domain $\Omega\subseteq\mathbb{R}^3$.
For convenience, the notation $\langle\cdot,\cdot\rangle$ represent the standard inner product in $L^2$. $[\partial^m,f]g$ denotes $\partial^m(fg)-f\partial^m g$. The same letter $C \geq 1$ represents a generic positive constant that depends on some known constants but is independent of $\epsilon$. $A\lesssim B$ means $A\leq CB$ for some constant $C>0$. Moreover, $A\approx B$ means that $c_1A\leq B\leq c_2A$ for some constants $c_1,c_2>0$, independent of $\epsilon$.

In addition, we introduce the following energy functional
\begin{align*}
E(\phi^{\epsilon},u^{\epsilon},\eta^{\epsilon},\tau^{\epsilon})
=\|\phi^{\epsilon}\|_{H^3_{P'(\rho^{\epsilon})}}^2
+\|u^{\epsilon}\|_{H^3_{\rho^{\epsilon}}}^2
+[\beta(L-1)+2\bar{\mathfrak{z}}]\|(\eta^{\epsilon}-1)\|_{H^3}^2+\frac{\beta}{2k^2}\|\tau^{\epsilon}\|_{H^3}^2,
\end{align*}
and the energy dissipative functional
\begin{align*}
D(u^{\epsilon},\eta^{\epsilon},\tau^{\epsilon})
=\mu_1\|\nabla u^{\epsilon}\|_{H^3}^2+\mu_2\|\mathrm{div} u^{\epsilon}\|_{H^3}^2
+\nu[\beta(L-1)+2\bar{\mathfrak{z}}]\|\nabla\eta^{\epsilon}\|_{H^3}^2+\frac{\beta A_0}{4k^2}\|\tau^{\epsilon}\|_{H^3}^2
+\frac{\beta\nu}{2k^2}\|\nabla\tau^{\epsilon}\|_{H^3}^2.
\end{align*}

\medskip
Now, we state the main results of this paper.
\begin{theorem}(uniform regularity)\label{Theorem3} Let $0<\epsilon\leq 1$, and consider the compressible model \eqref{ma1} with the initial data
\begin{align}
(\phi^{\epsilon}, u^{\epsilon}, \eta^{\epsilon},\tau^{\epsilon})\big|_{t=0}=(\phi_0^{\epsilon}, u_0^{\epsilon}, \eta_0^{\epsilon},\tau_0^{\epsilon}).\label{initial}
\end{align}
Assume that the initial data  $(\phi^{\epsilon}_0, u^{\epsilon}_0, \eta^{\epsilon}_0, \tau^{\epsilon}_0)$ satisfy $\rho^{\epsilon}_0=1+\epsilon\phi^{\epsilon}_0$, $\|\phi^{\epsilon}_0\|_{L^\infty}\leq\frac{1}{2}$, and
\begin{align}
 \|\phi^{\epsilon}_0\|_{H^3_{P'(\rho^{\epsilon}_0)}}^2 +\|u^{\epsilon}_0\|_{H^3_{\rho^{\epsilon}_0}}^2
+\|\eta^{\epsilon}_0-1\|_{H^3}^2
+\|\tau^{\epsilon}_0\|_{H^3}^2\leq\delta,\label{initial2}
\end{align}
for some small constant $\delta>0$ that does not depend on $\epsilon$. Then the system \eqref{ma1} and \eqref{initial} admits a unique global solution $(\phi^{\epsilon}, u^{\epsilon}, \eta^{\epsilon}, \tau^{\epsilon})$ which satisfies the following uniform (in $\epsilon$) bounds:
\begin{align}
&\|\phi^{\epsilon}\|_{ L^{\infty}(\mathbb{R}^+;H^3_{P'(\rho^{\epsilon})})}^2
+\|u^{\epsilon}\|_{ L^{\infty}(\mathbb{R}^+;H^3_{\rho^{\epsilon}})}^2
+\|\eta^{\epsilon}-1\|_{ L^{\infty}(\mathbb{R}^+;H^3)}^2+\|\tau^{\epsilon}\|_{ L^{\infty}(\mathbb{R}^+;H^3)}^2\label{clc-02}\\
&+\|\nabla u^{\epsilon}\|_{ L^2(\mathbb{R}^+;H^3)}^2
+\|\nabla\eta^{\epsilon}\|_{ L^2(\mathbb{R}^+;H^3)}^2
+\|\tau^{\epsilon}\|_{ L^2(\mathbb{R}^+;H^3)}^2
+\|\nabla\tau^{\epsilon}\|_{ L^2(\mathbb{R}^+;H^3)}^2
\lesssim \delta,\notag
\end{align}
\begin{align}
\|\rho^{\epsilon}\|_{L^{\infty}(\mathbb{R}^+;L^{\infty})}=
\|1+\epsilon\phi^{\epsilon}\|_{L^{\infty}(\mathbb{R}^+;L^{\infty})}
\approx 1,\label{clc-03}
\end{align}
and
\begin{align}
\|\partial_t\eta^{\epsilon}\|_{L^{\infty}(\mathbb{R}^+;H^1)}+
\|\partial_t\tau^{\epsilon}\|_{L^{\infty}(\mathbb{R}^+;H^1)}
\leq C .\label{clc-04}
\end{align}
Furthermore, suppose that
\begin{align}
\|\mathrm{div}u^{\epsilon}_0\|_{H^1}\leq C\epsilon,\;\;\;\|\nabla \phi^{\epsilon}_0\|_{H^1}\leq C\epsilon.\label{clc-05}
\end{align}
Then, for any $T>0$, there exists a positive constant $C(T)$ independent of $\epsilon$, such that
\begin{align}
\|\partial_t\phi^{\epsilon}\|_{L^{\infty}(0,T;H^1_{P'(\rho^{\epsilon})})}+
\|\partial_t u^{\epsilon}\|_{L^{\infty}(0,T;H^1_{\rho^{\epsilon}})}
\leq C(T) ,\label{clc-06}
\end{align}
and
\begin{align}
\frac{1}{\epsilon}\|\mathrm{div}u^{\epsilon}\|_{L^{\infty}(0,T;H^1)}+
\frac{1}{\epsilon}\|P'(\rho^{\epsilon})\nabla\phi^{\epsilon} \|_{L^{\infty}(0,T;H^1)}
\leq C(T).\label{clc-07}
\end{align}
\end{theorem}


The last result is about the limit from the global solution of the compressible Oldroyd-B model to the global solution of the incompressible counterpart and the convergence rate of the limit.
\begin{theorem}\label{Theorem2}(convergence rates) Consider the strong solutions $(\rho^{\epsilon}, u^{\epsilon}, \eta^{\epsilon}, \tau^{\epsilon})$ of the system \eqref{ma1} with the initial data \eqref{initial} constructed in Theorem \ref{Theorem3}, and the incompressible model \eqref{limit1} with the initial data
\begin{align}
( u, \eta,\tau)\big|_{t=0}=( u_0, \eta_0,\tau_0).\label{limit-initial}
\end{align}
Suppose in addition that the initial data satisfies $u_0,\eta_0-1,\tau_0\in H ^3$, $\mathrm{div}u_0=0$ and
\begin{align}
(u^{\epsilon}_0,\eta^{\epsilon}_0-1,\tau^{\epsilon}_0)\rightarrow (u_0,\eta_0-1,\tau_0) \;\;\text{strongly in}\;\; H^3, \;\;\text{as}\;\; \epsilon\rightarrow0.\label{clc-08}
\end{align}
Then there exists a subsequence $(\rho^{\epsilon}, u^{\epsilon}, \eta^{\epsilon}, \tau^{\epsilon})$ such that for any $T>0$,
\begin{align*}
&\rho^{\epsilon}\rightarrow1\;\;\text{strongly in}\;\; L^\infty(\mathbb{R}^+;H^3),\\
&\frac{1}{\epsilon^2}\nabla P(\rho^{\epsilon})
+\nabla\big(\beta(L-1)\eta^{\epsilon}+\bar{\mathfrak{z}}(\eta^{\epsilon})^{2}\big)
\rightarrow \nabla\pi\;\; \text{weakly-}*\;\text{in}\;\;L^\infty(0,T;H^1),\\
&(u^{\epsilon},\eta^{\epsilon}-1,\tau^{\epsilon})\rightarrow (u,\eta-1,\tau) \;\; \text{weakly-}*\;\text{in}\;\;L^\infty(\mathbb{R}^+;H^3)
\;\; \text{and}\;\;\text{strongly in}\;\; \mathcal{C}([0,T];H^2_{\text{loc}}),
\end{align*}
as $\epsilon\rightarrow0$, where $(u,\pi,\eta,\tau)$ is the global solution to the incompressible Oldroyd-B model \eqref{limit1} with the initial data \eqref{limit-initial}. Moreover, the following global energy bound holds:
\begin{align}
  &\|u\|_{ L^{\infty}(\mathbb{R}^+;H^3)}^2
  +\|\eta-1\|_{ L^{\infty}(\mathbb{R}^+;H^3)}^2+\|\tau\|_{ L^{\infty}(\mathbb{R}^+;H^3)}^2
  \label{clc-09}\\
  &+\|\nabla u\|_{ L^2(\mathbb{R}^+;H^3)}^2
  +\|\nabla\eta\|_{ L^2(\mathbb{R}^+;H^3)}^2
  +\|\tau\|_{ L^2(\mathbb{R}^+;H^3)}^2
  +\|\nabla\tau\|_{ L^2(\mathbb{R}^+;H^3)}^2
  \lesssim \delta,\notag
\end{align}
where $\delta$ are given in \eqref{initial2}.

Moreover, assume that
\begin{align}
\|\sqrt{\rho^{\epsilon}_0}u^{\epsilon}_0-u_0\|_{L^2}^2
+\|\eta^{\epsilon}_0-\eta_0\|_{L^2}^2+\|\tau^{\epsilon}_0-\tau_0\|_{L^2}^2
+\langle\Pi^{\epsilon}_0,1\rangle\lesssim\epsilon^{\alpha_0},\label{clc-10}
\end{align}
for some constant $\alpha_0>0$ independent of $\epsilon$, where $\Pi^{\epsilon}_0=\frac{1}{\epsilon^2}\frac{a}{\gamma-1}[(\rho^{\epsilon}_0)^{\gamma}-\gamma(\rho^{\epsilon}_0-1)-1]$. Then, for any fixed $T>0$, we have
\begin{align}
\|\sqrt{\rho^{\epsilon}}u^{\epsilon}-u\|_{L^2}^2
+\|\eta^{\epsilon}-\eta\|_{L^2}^2+\|\tau^{\epsilon}-\tau\|_{L^2}^2
+\langle\Pi^{\epsilon},1\rangle\leq C_T\epsilon^{\beta_0},\label{clc-11}
\end{align}
for all $t\in[0,T]$, where $\Pi^{\epsilon}=\frac{1}{\epsilon^2}\frac{a}{\gamma-1}[(\rho^{\epsilon})^{\gamma}-\gamma(\rho^{\epsilon}-1)-1]$, the constants $\beta_0=\min\{2,\alpha_0,1+\frac{\alpha_0}{2}\}>0$ and $C_T=C(1+T)\exp(CT)>0$ for some positive constant $C$ independent of $\epsilon$.
\end{theorem}

At the end of this subsection, we would like to sketch the main idea of this work. Inspired by \cite{FZ,HDW,LZ}, the local existence of strong solutions can first be obtained by using an iteration argument. Then, we seek some {\it a priori} estimates independent of $\epsilon$, so that the global existence of strong solution can be obtained by using the standard continuity arguments. To close the {\it a priori} assumption, we introduce some weighted norms to deal with the singular terms $\frac{1}{\epsilon}\mathrm{div}u^{\epsilon}$ and $\frac{1}{\epsilon}\frac{P'(\rho^{\epsilon})}{\rho^{\epsilon}}\nabla\phi^{\epsilon}$. When it comes to the incompressible limit of the solution, it is necessary to obtain the uniform estimate for their time derivatives. While justifying Theorem \ref{Theorem2}, the compactness arguments depending on Aubin-Lions-Simon Theorem (see Lemma \ref{aubin}) and the relative entropy method play a crucial role.

\medskip
In the last subsection, for reader's convenience, we present some known results that will be used later.

\subsection{Some known results}

\begin{lemma}(Moser-type inequality, \cite{MB}). For functions $f,g\in H^m\cap L^\infty$, $m\in \mathbb{Z}_+\cup\{0\}$, we have
\begin{align}
\|\nabla^m(fg)\|_{L^2}&\lesssim \|f\|_{L^\infty}\|\nabla^m g\|_{L^2}+\|\nabla^m f\|_{L^2}\|g\|_{L^\infty},\label{m1}\\
\sum_{0\leq|\alpha|\leq m}\|[\nabla^\alpha,f]g\|_{L^2}&\lesssim\|\nabla f\|_{L^\infty}\|\nabla^{m-1} g\|_{L^2}+\|\nabla^m f\|_{L^2}\|g\|_{L^\infty}.\label{m2}
\end{align}
In particular, if $m>\frac{3}{2}$, then
\begin{align*}
&\|fg\|_{H^m}\lesssim \|f\|_{H^m}\|g\|_{H^m}.
\end{align*}
\end{lemma}

\begin{lemma}\label{aubin}(Aubin-Lions-Simon Theorem, \cite{BF,S,S1}). Let $B_0\subset B_1\subset B_2$ be three Banach spaces. We assume that the embedding of $B_1$ in $B_2$ is continuous and that the embedding of $B_0$ in $B_1$ is compact. Let $1\leq p,r\leq +\infty$. For $T>0$, we define
\[E_{p,r}=\{u\in L^p(0,T; B_0), \partial_t u\in L^r(0,T;B_2)\}.\]
Then we have the following results:
\begin{itemize}
  \item [(1)] If $p<+\infty$, the embedding of $E_{p,r}$ in $L^p(0,T; B_1)$ is compact.
  \item [(2)] If $p=+\infty$ and $r>1$, the embedding of $E_{p,r}$ in $\mathcal{C}(0,T; B_1)$ is compact.
\end{itemize}
\end{lemma}

The rest of the paper is organized as follows. In Section \ref{global}, we obtain the global existence of system \eqref{ma1} with small initial data independent of $\epsilon$ by using the standard continuity arguments, as well as some uniform estimates of the solution. In Section \ref{sec4}, based on Theorem \ref{Theorem3}, we prove that the global strong solution of \eqref{ma1} converges to the global solution of the limiting target incompressible system \eqref{limit1} by compactness arguments. In addition, we also study the convergence rate in this section.

\section{Proof of Theorem \ref{Theorem3}}\label{global}
In this section, the global existence of system \eqref{ma1} and \eqref{initial} will be established by combining the local existence result with some global a priori estimates, with the help of the standard continuity arguments. At the same time, we are also going to get the uniform estimates \eqref{clc-02} and \eqref{clc-03}. Then, based on \eqref{clc-02}, \eqref{clc-03} and the additional condition \eqref{clc-05}, the uniform estimates of \eqref{clc-04}, \eqref{clc-06} and \eqref{clc-07} will be achieved step by step.

Actually, given some suitable initial data, the local existence and uniqueness of the strong solutions to \eqref{ma1} and \eqref{initial} can be obtained by using a standard iteration argument $($see $\cite{FZ,HDW,LZ}$ for instance$)$. We state the specific results bellow and omit the detail of the proof for brevity.
\begin{proposition}[local existence and uniqueness]\label{proposition1} Assume that $(\phi^{\epsilon}_0,u^{\epsilon}_0,\eta^{\epsilon}_0-1,\tau^{\epsilon}_0)\in H^3$. Then, there exists $T_0 > 0$ depending on $\epsilon$ and $\|(\phi^{\epsilon}_0, u^{\epsilon}_0, \eta^{\epsilon}_0-1, \tau^{\epsilon}_0)\|_{H^3}$ such that the initial-value problem \eqref{ma1} and \eqref{initial} has a unique strong solution $(\phi^{\epsilon}, u^{\epsilon}, \eta^{\epsilon}, \tau^{\epsilon})$ over $\mathbb{R}^3\times [0,T_0]$, which satisfies
\begin{equation*}
\begin{split}
&\phi^{\epsilon}\in
\mathcal{C}([0,T_0];H^3),\ \phi^{\epsilon}_t\in \mathcal{C}([0,T_0]; H^{2}),\\
& (u^{\epsilon},\eta^{\epsilon}-1,\tau^{\epsilon})\in \mathcal{C}([0,T_0];H^3)\cap L^2(0,T_0;H^{4}),\\
&(u^{\epsilon}_t,\eta^{\epsilon}_t, \tau^{\epsilon}_t)\in  \mathcal{C}([0,T_0];H^1)\cap L^2(0,T_0;H^{2}).
\end{split}
\end{equation*}
\end{proposition}

\medskip

To conclude the local solution is indeed a global one, it is crucial to prove the following proposition:
\begin{proposition}\label{4-proposition5}({\it A\; priori\;} estimate) Suppose that   $(\phi^{\epsilon},u^{\epsilon},\eta^{\epsilon},\tau^{\epsilon})$ is the strong solution of the Cauchy problem \eqref{ma1} and \eqref{initial} in $[0,T]$, where $T$ is a positive constant. Under the assumptions of Theorem \ref{Theorem3},  there exists a small positive
constant $\delta_1$ independent of $\epsilon$ and $T$, such that if
\begin{align}
E(\phi^{\epsilon},u^{\epsilon},\eta^{\epsilon},\tau^{\epsilon})\leq\delta_1,\label{s78-}
\end{align}
for any $t\in[0,T]$, then it holds that
\begin{align}
E(\phi^{\epsilon},u^{\epsilon},\eta^{\epsilon},\tau^{\epsilon})\leq\frac{\delta_1}{2}\label{s79},
\end{align}
for any $t\in[0,T]$.
\end{proposition}
\begin{remark}
The number $\delta_1$ is at least bigger than $4\delta$ $($determined by \eqref{4-s35}$)$.
\end{remark}

Proposition \ref{4-proposition5} can be derived from Lemmas \ref{4-Proposition0}-\ref{4-Proposition2} below.
\begin{lemma}\label{4-Proposition0} Under the same assumptions of Proposition \ref{4-proposition5}, we have
\begin{align}\label{4-s19-}
\frac{1}{2}\frac{d}{dt}E(\phi^{\epsilon}, u ^{\epsilon},\eta^{\epsilon},\tau^{\epsilon})
+\frac{3}{4}D( u ^{\epsilon},\eta^{\epsilon},\tau^{\epsilon})
\lesssim  \delta_1\|\nabla\phi^{\epsilon}\|_{H^2}^2.
\end{align}
\end{lemma}
\pf
Applying derivatives $\nabla^\ell(\ell=0,1,2,3)$ to the system \eqref{ma1}, taking inner product with $P'(\rho^{\epsilon})\nabla^\ell\phi^{\epsilon}$, $\rho^{\epsilon}\nabla^\ell  u ^{\epsilon}$, $[\beta(L-1)+2\bar{\mathfrak{z}}]\nabla^\ell (\eta^{\epsilon}-1)$ and $\frac{\beta}{2k^2}\nabla^\ell \tau^{\epsilon}$ respectively, and then adding the results, we can obtain
\begin{align}
\frac{1}{2}\frac{d}{dt}E(t)+D(t)
=&-\frac{1}{\epsilon}\sum\limits_{\ell=0}^{3}\big\langle\nabla^\ell
\mathrm{div} u ^{\epsilon},P'(\rho^{\epsilon})\nabla^\ell\phi^{\epsilon}\big\rangle
-\frac{1}{\epsilon}\sum\limits_{\ell=0}^{3}\big\langle\nabla^\ell\big(\frac{P'(\rho^{\epsilon})}{\rho^{\epsilon}}\nabla\phi^{\epsilon}\big),
\rho^{\epsilon}\nabla^{\ell} u ^{\epsilon}\big\rangle\notag\\
&+\frac{1}{2}\sum\limits_{\ell=0}^{3}\big\langle \partial_t P'(\rho^{\epsilon}), |\nabla^\ell\phi^{\epsilon}|^2\big\rangle
+\frac{1}{2}\sum\limits_{\ell=0}^{3}\big\langle \partial_t \rho^{\epsilon}, |\nabla^\ell  u ^{\epsilon}|^2\big\rangle\label{4-a}\\
&-\sum\limits_{\ell=0}^{3}\big\langle\nabla^\ell( u ^{\epsilon}\cdot\nabla\phi^{\epsilon}
+\phi^{\epsilon}\mathrm{div} u ^{\epsilon}),P'(\rho^{\epsilon})\nabla^\ell\phi^{\epsilon}\big\rangle
-\sum\limits_{\ell=0}^{3}\big\langle\nabla^\ell( u ^{\epsilon}\cdot\nabla  u ^{\epsilon})
,\rho^{\epsilon}\nabla^\ell  u ^{\epsilon}\big\rangle\notag\\
&-\Big(\beta(L-1) +2\bar{\mathfrak{z}}\Big)\sum\limits_{\ell=0}^{3}\big\langle\nabla^\ell\left( u ^{\epsilon}\cdot\nabla \eta^{\epsilon}
+(\eta^{\epsilon}-1)\mathrm{div} u ^{\epsilon}\right),\nabla^\ell(\eta^{\epsilon}-1)\big\rangle\notag\\
&-\frac{\beta}{2k^2}\sum\limits_{\ell=0}^{3}\big\langle\nabla^\ell( u ^{\epsilon}\cdot\nabla\tau^{\epsilon}
+\tau^{\epsilon}\mathrm{div} u ^{\epsilon})-\nabla^\ell(\nabla  u ^{\epsilon}\tau^{\epsilon}
+\tau^{\epsilon}\nabla^T  u ^{\epsilon}),\nabla^\ell\tau^{\epsilon}\big\rangle\notag
  \\
&+\frac{\beta}{2k}\sum\limits_{\ell=0}^{3}\big\langle \nabla^\ell\big((\eta^{\epsilon}-1)(\nabla  u ^{\epsilon}
+\nabla^T  u ^{\epsilon})\big),\nabla^\ell\tau^{\epsilon}\big\rangle
-2\bar{\mathfrak{z}}\sum\limits_{\ell=0}^{3}\big\langle \nabla^{\ell}\big((\eta^{\epsilon}-1)
\nabla\eta^{\epsilon}\big),\nabla^{\ell} u ^{\epsilon}\big\rangle\notag\\
 &-\sum\limits_{\ell=1}^{3}\big\langle[\nabla^{\ell}, (\frac{1}{\rho^{\epsilon}}-1)]\big(\beta(L-1)+2\bar{\mathfrak{z}}\eta^{\epsilon}\big)\nabla\eta^{\epsilon}
,\rho^{\epsilon}\nabla^{\ell} u ^{\epsilon}\big\rangle\notag\\
&+\mu_1\sum\limits_{\ell=1}^{3}\big\langle[\nabla^{\ell}, (\frac{1}{\rho^{\epsilon}}-1)]\Delta u ^{\epsilon},\rho^{\epsilon}\nabla^{\ell} u ^{\epsilon}\big\rangle
+\mu_2 \sum\limits_{\ell=1}^{3}\big\langle[\nabla^{\ell}, (\frac{1 }{\rho^{\epsilon}}-1)]\nabla\mathrm{div} u ^{\epsilon},\rho^{\epsilon}\nabla^{\ell} u ^{\epsilon}\big\rangle\notag\\
&+\frac{\beta}{k}\sum\limits_{\ell=1}^{3}\big\langle[\nabla^{\ell}, (\frac{1}{\rho^{\epsilon}}-1)]\mathrm{div}\tau^{\epsilon}
,\rho^{\epsilon}\nabla^{\ell} u ^{\epsilon}\big\rangle
:=\sum\limits_{j=1}^{14}I_{i}.\notag
\end{align}
Next, we turn to deal with the terms $I_{1}$-$I_{14}$. Firstly, for $I_1$ and $I_2$, thanks to integration by parts, we have
\begin{align*}
I_1+I_2=&\frac{1}{\epsilon}\sum\limits_{\ell=0}^{3}\big\langle\nabla^\ell
 u ^{\epsilon},\nabla\big(P'(\rho^{\epsilon})\nabla^\ell\phi^{\epsilon}\big)\big\rangle
-\frac{1}{\epsilon}\sum\limits_{\ell=0}^{3}\big\langle P'(\rho^{\epsilon})\nabla^{\ell+1}\phi^{\epsilon},
\nabla^{\ell} u ^{\epsilon}\big\rangle\\
&-\frac{1}{\epsilon}\sum\limits_{\ell=0}^{3}\big\langle[\nabla^\ell,(\frac{P'(\rho^{\epsilon})}{\rho^{\epsilon}}-a\gamma)]\nabla\phi^{\epsilon},
\rho^{\epsilon}\nabla^{\ell} u ^{\epsilon}\big\rangle\\
=&\sum\limits_{\ell=0}^{3}\big\langle\nabla^\ell
 u ^{\epsilon},P''(\rho^{\epsilon})\nabla\phi^{\epsilon}\nabla^\ell\phi^{\epsilon}\big\rangle
-\frac{1}{\epsilon}\sum\limits_{\ell=0}^{3}\big\langle[\nabla^\ell,(\frac{P'(\rho^{\epsilon})}{\rho^{\epsilon}}-a\gamma)]\nabla\phi^{\epsilon},
\rho^{\epsilon}\nabla^{\ell} u ^{\epsilon}\big\rangle.
\end{align*}
Further, H\"older inequality, Sobolev inequality, Cauchy inequality, \eqref{m2} and \eqref{s78-} imply the following estimate,
\begin{align}
|I_1+I_2|\lesssim& \|P''(\rho^{\epsilon})\|_{L^\infty}\|\nabla\phi^{\epsilon}\|_{H^2}
\| u ^{\epsilon}\|_{H^3}\|\nabla\phi^{\epsilon}\|_{H^2}\label{4-a1}\\
&+\frac{1}{\epsilon}\big\|\nabla\frac{P'(\rho^{\epsilon})}{\rho^{\epsilon}}\big\|_{H^2}
\|\nabla\phi^{\epsilon}\|_{H^2}\|\rho^{\epsilon}\|_{L^\infty}\|\nabla u ^{\epsilon}\|_{H^2}\notag\\
\lesssim& \delta_1\|\nabla\phi^{\epsilon}\|_{H^2}^2+\delta_1\|\nabla u ^{\epsilon}\|_{H^2}^2,\notag
\end{align}
where we have used the following technique to cancel the singularity of $\frac{1}{\epsilon}$. Based on $\rho^{\epsilon}=1+\epsilon\phi^{\epsilon}$, H\"older inequality, Sobolev inequality and \eqref{s78-},
\begin{align*}
\frac{1}{\epsilon}\big\|\nabla\frac{P'(\rho^{\epsilon})}{\rho^{\epsilon}}\big\|_{H^2}
=&\frac{1}{\epsilon}\big\| \frac{P''(\rho^{\epsilon})\epsilon\nabla\phi^{\epsilon}}{\rho^{\epsilon}}
-\frac{P'(\rho^{\epsilon})\epsilon\nabla\phi^{\epsilon}}{(\rho^{\epsilon})^2}\big\|_{H^2}\\
=& \big\| \frac{P''(\rho^{\epsilon}) \nabla\phi^{\epsilon}}{\rho^{\epsilon}}
-\frac{P'(\rho^{\epsilon}) \nabla\phi^{\epsilon}}{(\rho^{\epsilon})^2}\big\|_{H^2}\\
\lesssim&\delta_1\|\nabla\phi^{\epsilon}\|_{H^2}.
\end{align*}
For $I_3$, we recall that
$\partial_t\rho^{\epsilon}=- u ^{\epsilon}\cdot\nabla\rho^{\epsilon}
-\rho^{\epsilon}\mathrm{div} u ^{\epsilon}$, and then multiplying the equality by $P''(\rho^{\epsilon})$ to obtain
\begin{align}\label{4-s14}
\partial_tP'(\rho^{\epsilon})=- u ^{\epsilon}\cdot\nabla P'(\rho^{\epsilon})
-P''(\rho^{\epsilon})\rho^{\epsilon}\mathrm{div} u ^{\epsilon}.
\end{align}
Using H\"older inequality, Sobolev inequality, Cauchy inequality and \eqref{s78-}, we can estimate $I_3$ as follows
\begin{align}\label{4-a2}
I_3=&\frac{1}{2}\sum\limits_{\ell=0}^{3}\big\langle - u ^{\epsilon}\cdot\nabla P'(\rho^{\epsilon})
-P''(\rho^{\epsilon})\rho^{\epsilon}\mathrm{div} u ^{\epsilon}
, |\nabla^\ell\phi^{\epsilon}|^2\big\rangle\notag\\
\lesssim&\|P''(\rho^{\epsilon})\|_{L^\infty}(\| u ^{\epsilon}\|_{L^\infty}
\|\epsilon\nabla\phi^{\epsilon}\|_{L^2}+\|\rho^{\epsilon}\|_{L^\infty}
\|\mathrm{div} u ^{\epsilon}\|_{L^2})
\|\phi^{\epsilon}\|_{L^3}\|\phi^{\epsilon}\|_{L^6}\notag\\
&+\|P''(\rho^{\epsilon})\|_{L^\infty}(\| u ^{\epsilon}\|_{L^\infty}
\|\epsilon\nabla\phi^{\epsilon}\|_{L^\infty}+\|\rho^{\epsilon}\|_{L^\infty}
\|\mathrm{div} u ^{\epsilon}\|_{L^\infty})
\|\nabla\phi^{\epsilon}\|_{H^2}^2\\
\lesssim&\delta_1\|\nabla u ^{\epsilon}\|_{H^2}^2+\delta_1
\|\nabla\phi^{\epsilon}\|_{H^2}^2\notag.
\end{align}
For $I_4$ and $I_6$, thanks to integration by parts and $\partial_t\rho^{\epsilon}+\mathrm{div}(\rho^{\epsilon} u ^{\epsilon})=0$, combining with H\"older inequality, Sobolev inequality, Cauchy inequality and \eqref{s78-}, we have
\begin{align}\label{4-a3}
I_4+I_6=&\frac{1}{2}\sum\limits_{\ell=0}^{3}\big\langle \partial_t \rho^{\epsilon}, |\nabla^\ell  u ^{\epsilon}|^2\big\rangle
-\sum\limits_{\ell=0}^{3}\big\langle\nabla^\ell( u ^{\epsilon}\cdot\nabla  u ^{\epsilon})
,\rho^{\epsilon}\nabla^\ell  u ^{\epsilon}\big\rangle\notag\\
=&-\sum\limits_{\ell=0}^{3}\big\langle[\nabla^\ell, u ^{\epsilon}\cdot\nabla] u ^{\epsilon},\rho^{\epsilon}\nabla^\ell u ^{\epsilon}\big\rangle
\\
\lesssim&\|\nabla u ^{\epsilon}\|_{H^2}^2
\|\rho^{\epsilon}\|_{L^\infty}
\| u ^{\epsilon}\|_{H^3}\notag\\
\lesssim&\delta_1\|\nabla u ^{\epsilon}\|_{H^2}^2\notag.
\end{align}
When it comes to $I_5$, it seems impossible to estimate those terms containing $\nabla^4\phi^{\epsilon}$ directly since the equation of $\phi^{\epsilon}$ has no dissipative term. For this reason, we use integration by parts to transfer the derivative to other terms. More precisely, from H\"older inequality, Sobolev inequality, \eqref{m1}, \eqref{m2} and \eqref{s78-},  we obtain
\begin{align}
I_5=&-\sum\limits_{\ell=0}^{3}\big\langle[\nabla^\ell, u ^{\epsilon}\cdot\nabla]\phi^{\epsilon}+ u ^{\epsilon}\cdot\nabla^{\ell+1}\phi^{\epsilon}
+\nabla^\ell(\phi^{\epsilon}\mathrm{div} u ^{\epsilon}),P'(\rho^{\epsilon})\nabla^\ell\phi^{\epsilon}\big\rangle\notag\\
=&-\sum\limits_{\ell=0}^{3}\big\langle[\nabla^\ell, u ^{\epsilon}\cdot\nabla]\phi^{\epsilon}
+\nabla^\ell(\phi^{\epsilon}\mathrm{div} u ^{\epsilon}),P'(\rho^{\epsilon})\nabla^\ell\phi^{\epsilon}\big\rangle
+\frac{1}{2}\sum\limits_{\ell=0}^{3}\big\langle\mathrm{div}\big(P'(\rho^{\epsilon}) u ^{\epsilon}\big),|\nabla^\ell\phi^{\epsilon}|^2\big\rangle\notag\\
\lesssim&(\|\nabla  u ^{\epsilon}\|_{H^2}+\|\mathrm{div} u ^{\epsilon}\|_{H^3})\|\nabla \phi^{\epsilon}\|_{H^2}
\|P'(\rho^{\epsilon})\|_{L^\infty}\|\phi^{\epsilon}\|_{H^3}\label{4-a4}\\
&+\left(\|\mathrm{div} u ^{\epsilon}\|_{H^2}\|P'(\rho^{\epsilon})\|_{L^\infty}
+\| u ^{\epsilon}\|_{H^2}\|P''(\rho^{\epsilon})\|_{L^\infty}\|\epsilon\nabla\phi^{\epsilon}\|_{L^\infty}
\right)\|\phi^{\epsilon}\|_{H^3}\|\nabla\phi^{\epsilon}\|_{H^2}\notag\\
\lesssim&\delta_1\|\nabla u ^{\epsilon}\|_{H^3}^2+\delta_1
\|\nabla\phi^{\epsilon}\|_{H^2}^2\notag.
\end{align}
The remaining terms in \eqref{4-a} can be estimated by using H\"older inequality, Sobolev inequality, Cauchy inequality, \eqref{m1}, \eqref{m2} and \eqref{s78-} as follows
\begin{align}\label{4-a6}
&|I_7+I_8+I_9+I_{10}|\notag\\
\lesssim&(\| u ^{\epsilon}\|_{H^3}\|\nabla \eta^{\epsilon}\|_{H^3}+\|(\eta^{\epsilon}-1)\|_{H^3}\|\nabla  u ^{\epsilon}\|_{H^3})\|\nabla\eta^{\epsilon}\|_{H^2}\notag\\
&+\big[\| u ^{\epsilon}\|_{H^3}\|\nabla \tau^{\epsilon}\|_{H^3}
+\big(\|\tau^{\epsilon}\|_{H^3}+\|(\eta^{\epsilon}-1)\|_{H^3}\big)\|\nabla u ^{\epsilon}\|_{H^3}\big]\|\nabla\tau^{\epsilon}\|_{H^2}\\
&+\|(\eta^{\epsilon}-1)\|_{H^3}\|\nabla \eta^{\epsilon}\|_{H^3}\|\nabla u ^{\epsilon}\|_{H^2} \notag\\
\lesssim&\delta_1\|\nabla u ^{\epsilon}\|_{H^3}^2+\delta_1\|\nabla\eta^{\epsilon}\|_{H^3}^2+\delta_1\|\nabla\tau^{\epsilon}\|_{H^3}^2
,\notag
\end{align}
and
\begin{align}
&|I_{11}+I_{12}+I_{13}+I_{14}|\notag\\
\lesssim&\|\nabla\frac{1}{\rho^{\epsilon}}\|_{H^2}
\big(\|\nabla^2 u ^{\epsilon}\|_{H^2}+ \|\nabla\eta^{\epsilon}\|_{H^2} +\|\nabla\tau^{\epsilon}\|_{H^2}
\big)\|\rho^{\epsilon}\|_{L^\infty}
\|\nabla u ^{\epsilon}\|_{H^2}\label{4-a7}\\
\lesssim&\delta_1\|\nabla u ^{\epsilon}\|_{H^3}^2+\delta_1\|\nabla\eta^{\epsilon}\|_{H^3}^2+\delta_1\|\nabla\tau^{\epsilon}\|_{H^3}^2
.\notag
\end{align}
Finally, putting \eqref{4-a1} and \eqref{4-a2}-\eqref{4-a7} into \eqref{4-a}, we can get directly
\begin{align*}
\frac{1}{2}\frac{d}{dt}E(t)+D(t)\lesssim\delta_1\|\nabla u ^{\epsilon}\|_{H^3}^2+\delta_1\|\nabla\eta^{\epsilon}\|_{H^3}^2+\delta_1\|\nabla\tau^{\epsilon}\|_{H^3}^2
+\delta_1\|\nabla\phi^{\epsilon}\|_{H^2}^2.
\end{align*}
Then, choosing $\delta_1$ sufficiently small, we can obtain \eqref{4-s19-}. The proof of Lemma \ref{4-Proposition0} is completed.
\endpf

\medskip
In the following lemma, we obtain some dissipation estimates of~$\phi^{\epsilon}$, which are independent of $\epsilon$.
\begin{lemma}\label{4-Proposition2} Under the same assumptions of Proposition \ref{4-proposition5}, we have
\begin{align}\label{4-2}
\sum\limits_{\ell=0}^{2}\frac{d}{dt}\langle\nabla^{\ell} u ^{\epsilon},\epsilon\nabla^{\ell+1}\phi^{\epsilon} \rangle
+\frac{3}{4}\|\nabla\phi^{\epsilon}\|_{H^2_{\frac{P'(\rho^{\epsilon})}{\rho^{\epsilon}}}}^2
\lesssim  D( u ^{\epsilon},\eta^{\epsilon},\tau^{\epsilon})+\|\mathrm{div} u ^{\epsilon}\|_{H^2}^2.
\end{align}
\end{lemma}
\pf
Let~$\ell=0,1,2$. Applying the derivative operator $\nabla^{\ell}$ to the equations \eqref{ma1}$_1$ and \eqref{ma1}$_2$, taking inner product with $-\epsilon \mathrm{div}\nabla^{\ell} u ^{\epsilon}$ and $\epsilon\nabla^{\ell+1}\phi^{\epsilon}$ respectively, and then adding
the results together, we can obtain
\begin{align}
&\sum\limits_{\ell=0}^{2}\frac{d}{dt}\langle\nabla^{\ell} u ^{\epsilon},\epsilon\nabla^{\ell+1}\phi^{\epsilon} \rangle
+\|\nabla\phi^{\epsilon}\|_{H^2_{\frac{P'(\rho^{\epsilon})}{\rho^{\epsilon}}}}^2
-\|\mathrm{div} u ^{\epsilon}\|_{H^2}^2\label{4-s18}\\
=&\sum\limits_{\ell=0}^{2}\langle\epsilon\mathrm{div}\nabla^{\ell} u ^{\epsilon},\nabla^{\ell}( u ^{\epsilon}\cdot\nabla\phi^{\epsilon}
+\phi^{\epsilon}\mathrm{div} u ^{\epsilon})\rangle
-\sum\limits_{\ell=0}^{2}\langle[\nabla^{\ell},\frac{P'(\rho^{\epsilon})}{\rho^{\epsilon}}]\nabla
\phi^{\epsilon}+\epsilon\nabla^{\ell}( u ^{\epsilon}\cdot \nabla u ^{\epsilon}),\nabla^{\ell+1}\phi^{\epsilon} \rangle\notag\\
&+\sum\limits_{\ell=0}^{2}\langle\nabla^{\ell}\left\{\frac{1}{\rho^{\epsilon}}
\big[-\nabla\big(\beta(L-1)\eta^{\epsilon}+\bar{\mathfrak{z}}(\eta^{\epsilon})^{2}\big)+\mu_1\Delta u ^{\epsilon}+\mu_2\nabla\mathrm{div} u  ^{\epsilon}+\frac{\beta}{k}\mathrm{div}\tau ^{\epsilon}\big]
\right\},\epsilon\nabla^{\ell+1}\phi^{\epsilon}\rangle.\notag
\end{align}
Based on H\"older inequality, Sobolev inequality and Cauchy inequality, the last line on the right-hand side of \eqref{4-s18} can be controlled as follows
\begin{align}\label{4-s20-}
&\sum\limits_{\ell=0}^{2}\langle\nabla^{\ell}\left[\frac{1}{\rho^{\epsilon}}
\bigg(-\nabla\big(\beta(L-1)\eta^{\epsilon}+\bar{\mathfrak{z}}(\eta^{\epsilon})^{2}\big)+\mu_1\Delta u ^{\epsilon}+\mu_2\nabla\mathrm{div} u  ^{\epsilon}+\frac{\beta}{k}\mathrm{div}\tau ^{\epsilon}\bigg)
\right],\epsilon\nabla^{\ell+1}\phi^{\epsilon}\rangle\notag\\
\lesssim& \Big(\| \frac{1}{\rho^{\epsilon}}\|_{L^\infty}+ \|\nabla\frac{1}{\rho^{\epsilon}}\|_{H^2}\Big)
\Big(
\big(1+\|(\eta^{\epsilon}-1)\|_{H^2}\big)\|\nabla\eta^{\epsilon}\|_{H^2}+\|\nabla^2 u ^{\epsilon}\|_{H^2}+\|\nabla\tau^{\epsilon}\|_{H^2}\Big)
\|\nabla\phi^{\epsilon}\|_{H^2} \notag\\
\lesssim&(1+\|\nabla\phi^{\epsilon}\|_{H^2}^3)
(\|\nabla\eta^{\epsilon}\|_{H^2}
+\|(\eta^{\epsilon}-1)\|_{H^2}\|\nabla\eta^{\epsilon}\|_{H^2}+\|\nabla^2 u ^{\epsilon}\|_{H^2}+\|\nabla\tau^{\epsilon}\|_{H^2})
\|\nabla\phi^{\epsilon}\|_{H^2}\notag\\
\leq& C_{\theta_1}\Big(\|\nabla u ^{\epsilon}\|_{H^3}^2+ \|\nabla\eta^{\epsilon}\|_{H^3}^2+ \|\nabla\tau^{\epsilon}\|_{H^3}^2\Big)
+\theta_1\|\nabla\phi^{\epsilon}\|_{H^2}^2,
\end{align}
where $\theta_1>0$ is a constant independent of $\epsilon$. Similarly, for the first term on the right-hand side of \eqref{4-s18}, we have
\begin{align}
&\sum\limits_{\ell=0}^{2}\langle\epsilon\mathrm{div}\nabla^{\ell} u ^{\epsilon},\nabla^{\ell}( u ^{\epsilon}\cdot\nabla\phi^{\epsilon}
+\phi^{\epsilon}\mathrm{div} u ^{\epsilon})\rangle\notag\\
&-\sum\limits_{\ell=0}^{2}\langle[\nabla^{\ell},(\frac{P'(\rho^{\epsilon})}{\rho^{\epsilon}}-a\gamma)]\nabla
\phi^{\epsilon}+\epsilon\nabla^{\ell}( u ^{\epsilon}\cdot \nabla u ^{\epsilon}),\nabla^{\ell+1}\phi^{\epsilon} \rangle\notag\\
\lesssim&(1+\|\phi^{\epsilon}\|_{H^3})
(\| u ^{\epsilon}\|_{H^3}+\|\phi^{\epsilon}\|_{H^3})
(\|\nabla  u ^{\epsilon}\|_{H^3}^2+\|\nabla\phi^{\epsilon}\|_{H^2}^2)\label{4-s20--}\\
\lesssim&\delta_1\|\nabla u ^{\epsilon}\|_{H^3}^2+\delta_1\|\nabla\phi^{\epsilon}\|_{H^2}^2.\notag
\end{align}
Hence, putting \eqref{4-s20-} and \eqref{4-s20--} into \eqref{4-s18}, we can get
\begin{align}\label{4-s20}
&\sum\limits_{\ell=0}^{2}\frac{d}{dt}\langle\nabla^{\ell} u ^{\epsilon},\epsilon\nabla^{\ell+1}\phi^{\epsilon} \rangle
+\|\nabla\phi^{\epsilon}\|_{H^2_{\frac{P'(\rho^{\epsilon})}{\rho^{\epsilon}}}}^2
-\|\mathrm{div} u ^{\epsilon}\|_{H^2}^2\\
\lesssim& (\delta_1+C_{\theta_1})\Big( \|\nabla u ^{\epsilon}\|_{H^3}^2+ \|\nabla\eta^{\epsilon}\|_{H^3}^2+ \|\nabla\tau^{\epsilon}\|_{H^3}^2\Big)
+(\delta_1 + \theta_1)\|\nabla\phi^{\epsilon}\|_{H^2}^2\notag\\
\lesssim&(\delta_1+C_{\theta_1}) D( u ^{\epsilon},\eta^{\epsilon},\tau^{\epsilon})
+(\delta_1 + \theta_1)\|\nabla\phi^{\epsilon}\|_{H^2_{\frac{P'(\rho^{\epsilon})}{\rho^{\epsilon}}}}^2,\notag
\end{align}
where we have used the fact that
\beq\label{4-3}
\|\nabla^\ell\phi^{\epsilon}\|_{L^2}\lesssim\big\|\frac{\rho^{\epsilon}}{P'(\rho^{\epsilon})}\big\|_{L^\infty}^\frac{1}{2}
\bigg\|\bigg(\frac{P'(\rho^{\epsilon})}{\rho^{\epsilon}}\bigg)^\frac{1}{2}\nabla^\ell\phi^{\epsilon}\bigg\|_{L^2}
\lesssim\|\nabla^\ell\phi^{\epsilon}\|_{L^2_{\frac{P'(\rho^{\epsilon})}{\rho^{\epsilon}}}}, \ \ell\in\mathbb{N}.
\eeq
Then, choosing $\delta_1$ and $\theta_1$ sufficiently small in \eqref{4-s20}, we can obtain
\begin{align*}
\sum\limits_{\ell=0}^{2}\frac{d}{dt}\langle\nabla^{\ell} u ^{\epsilon},\epsilon\nabla^{\ell+1}\phi^{\epsilon} \rangle
+\frac{3}{4}\|\nabla\phi^{\epsilon}\|_{H^2_{\frac{P'(\rho^{\epsilon})}{\rho^{\epsilon}}}}^2
\lesssim  D( u ^{\epsilon},\eta^{\epsilon},\tau^{\epsilon})+\|\mathrm{div} u ^{\epsilon}\|_{H^2}^2,
\end{align*}
which implies \eqref{4-2}. The proof of Lemma \ref{4-Proposition2} is completed.
\endpf

\medskip
Based on Lemmas \ref{4-Proposition0} and \ref{4-Proposition2}, we are ready to prove Proposition \ref{4-proposition5}.
\paragraph{Proof of Proposition \ref{4-proposition5}:}
Combining \eqref{4-s19-} with $\theta_2\times$\eqref{4-2}, where $\theta_2$ is a positive constant independent of $\epsilon$ and using \eqref{4-3}, we can deduce that
\begin{align}\label{4-s73}
&\frac{d}{dt}\bigg(\frac{1}{2}E(\phi^{\epsilon}, u ^{\epsilon},\eta^{\epsilon},\tau^{\epsilon})
+\theta_2\sum\limits_{\ell=0}^{2}\langle\nabla^{\ell} u ^{\epsilon},\epsilon\nabla^{\ell+1}\phi^{\epsilon} \rangle\bigg)
+\frac{3}{4}D( u ^{\epsilon},\eta^{\epsilon},\tau^{\epsilon})
+\frac{3\theta_2}{4}\|\nabla\phi^{\epsilon}\|_{H^2_{\frac{P'(\rho^{\epsilon})}{\rho^{\epsilon}}}}^2\notag\\
\lesssim&\delta_1\|\nabla\phi^{\epsilon}\|_{H^2}^2+\theta_2\|\mathrm{div} u ^{\epsilon}\|_{H^2}^2
+\theta_2  D( u ^{\epsilon},\eta^{\epsilon},\tau^{\epsilon})\\
\lesssim&\delta_1\|\nabla\phi^{\epsilon}\|_{H^2_{\frac{P'(\rho^{\epsilon})}{\rho^{\epsilon}}}}^2
+\theta_2\|\mathrm{div} u ^{\epsilon}\|_{H^2}^2
+\theta_2  D( u ^{\epsilon},\eta^{\epsilon},\tau^{\epsilon}).\notag
\end{align}
Choosing $\theta_2$ and $\delta_1$ sufficiently small in \eqref{4-s73}, then we obtain the following inequality
\begin{equation*}
\frac{d}{dt}\bigg(\frac{1}{2}E(\phi^{\epsilon}, u ^{\epsilon},\eta^{\epsilon},\tau^{\epsilon})
+\theta_2\sum\limits_{\ell=0}^{2}\langle\nabla^{\ell} u ^{\epsilon},\epsilon\nabla^{\ell+1}\phi^{\epsilon} \rangle\bigg)
+\frac{1}{2}D( u ^{\epsilon},\eta^{\epsilon},\tau^{\epsilon})
+\frac{\theta_2}{2}\|\nabla\phi^{\epsilon}\|_{H^2_{\frac{P'(\rho^{\epsilon})}{\rho^{\epsilon}}}}^2
\leq0.
\end{equation*}
Further, integrating the above inequality over $(0,t)$, we can get
\begin{align} \mathcal{E}(\phi^{\epsilon}, u ^{\epsilon},\eta^{\epsilon},\tau^{\epsilon})
+\frac{1}{2}\int_0^t\big(D( u ^{\epsilon},\eta^{\epsilon},\tau^{\epsilon})
+\theta_2\|\nabla\phi^{\epsilon}\|_{H^2_{\frac{P'(\rho^{\epsilon})}{\rho^{\epsilon}}}}^2\big)\,\mathrm{d}s
\leq\mathcal{E}(\phi^{\epsilon}, u ^{\epsilon},\eta^{\epsilon},\tau^{\epsilon})(0)\label{4-s78},
\end{align}
where
\[\mathcal{E}(\phi^{\epsilon}, u ^{\epsilon},\eta^{\epsilon},\tau^{\epsilon})
=\frac{1}{2}E(\phi^{\epsilon}, u ^{\epsilon},\eta^{\epsilon},\tau^{\epsilon})
+\theta_2\sum\limits_{\ell=0}^{2}\langle\nabla^{\ell} u ^{\epsilon},\epsilon\nabla^{\ell+1}\phi^{\epsilon} \rangle.\]
In addition, using the smallness of ~$\theta_2$, we are easy to check that $\mathcal{E}(\phi^{\epsilon}, u ^{\epsilon},\eta^{\epsilon},\tau^{\epsilon})$ is equivalent to $E(\phi^{\epsilon}, u ^{\epsilon},\eta^{\epsilon},\tau^{\epsilon})$. In other words, there exist positive constants $C_1$ and $C_2$ independent of $\epsilon$, such that
\begin{align*}
C_1E(\phi^{\epsilon}, u ^{\epsilon},\eta^{\epsilon},\tau^{\epsilon})
\leq \mathcal{E}(\phi^{\epsilon}, u ^{\epsilon},\eta^{\epsilon},\tau^{\epsilon})
\leq C_2E(\phi^{\epsilon}, u ^{\epsilon},\eta^{\epsilon},\tau^{\epsilon}),
\end{align*}
which, combining with \eqref{initial2} and \eqref{4-s78}, implies that
\begin{align*}
E(\phi^{\epsilon}, u ^{\epsilon},\eta^{\epsilon},\tau^{\epsilon})\leq \frac{C_2}{C_1}E(\phi^{\epsilon},u^{\epsilon},\eta^{\epsilon},\tau^{\epsilon})(0)
\leq\frac{C_2}{C_1}\delta.
\end{align*}
Let
\begin{align}\label{4-s35}
\delta\leq \frac{C_1}{2 C_2}\delta_1,
\end{align}
and then we can get~\eqref{s79}. The proof of Proposition \ref{4-proposition5} is finished.
\endpf
\begin{remark}
With the local existence result of Proposition \ref{proposition1} and the global {\it a priori} estimates of Proposition \ref{4-proposition5}, by using the standard continuity arguments, it is not hard to get the global existence of the solution. Moreover, we can also get the uniform estimate \eqref{4-s78} for all $t\in\mathbb{R}^+$, which means that the inequalities \eqref{clc-02} and \eqref{clc-03} are established.
\end{remark}

To prove Theorem \ref{Theorem3}, what left is to obtain the uniform bounds \eqref{clc-04}, \eqref{clc-06} and \eqref{clc-07}. This is achieved in the following proposition with the help of the uniform estimates \eqref{clc-02}, \eqref{clc-03} and the additional assumptions \eqref{clc-05}.
\begin{proposition}\label{4-proposition4} Under the same assumptions of Theorem \ref{Theorem3}, the solution $(\phi^{\epsilon}, u ^{\epsilon}, \eta^{\epsilon}, \tau^{\epsilon})$ to the
Cauchy problem \eqref{ma1} and \eqref{initial} admits the following uniform $(\text{in }\epsilon)$ bounds:
\begin{align}\label{4-4}
\|\partial_t\eta^{\epsilon}\|_{L^{\infty}(\mathbb{R}^+;H^1)}+
\|\partial_t\tau^{\epsilon}\|_{L^{\infty}(\mathbb{R}^+;H^1)}
\leq C,
\end{align}
and for any fixed $T>0$,  there exists a positive constant $C(T)$ independent of $\epsilon$, such that
\begin{align}
&\|\partial_t\phi^{\epsilon}\|_{L^{\infty}(0;T;H^1_{P'(\rho^{\epsilon})})}+
\|\partial_t  u ^{\epsilon}\|_{L^{\infty}(0;T;H^1_{\rho^{\epsilon}})}
\leq C(T) ,\label{4-5}\\
&\frac{1}{\epsilon}\|\mathrm{div} u ^{\epsilon}\|_{L^{\infty}(0;T;H^1)}+
\frac{1}{\epsilon}\|P'(\rho^{\epsilon})\nabla\phi^{\epsilon} \|_{L^{\infty}(0;T;H^1)}
\leq C(T).\label{4-6}
\end{align}
\end{proposition}
\pf
Firstly, for $\ell=0,1$, by virtue of the last two equations of \eqref{ma1}, the uniform estimates \eqref{clc-02}, \eqref{m1}, H\"older inequality and Sobolev inequality, it is not hard to get
\begin{align}
\|\partial_t\nabla^{\ell}\eta^{\epsilon}\|_{L^2}
&\lesssim\|\nabla^{\ell}\mathrm{div}(\eta ^{\epsilon}  u ^{\epsilon})\|_{L^2}+\|\Delta \nabla^{\ell}\eta^{\epsilon}\|_{L^2}\notag\\
&\lesssim (\|\nabla\eta ^{\epsilon}\|_{H^1}+\| \eta ^{\epsilon}\|_{L^\infty})\| u ^{\epsilon}\|_{H^2}+\|\nabla^2\eta ^{\epsilon}\|_{H^1}\label{4-s28}\\
&\lesssim \delta\notag
\end{align}
and
\begin{align}
\|\partial_t\nabla^{\ell}\tau^{\epsilon}\|_{L^2}
\lesssim&\|\nabla^{\ell}\mathrm{div} ( u ^{\epsilon}\tau^{\epsilon})\|_{L^2}+\|\nabla^{\ell}(\nabla  u ^{\epsilon}
\tau^{\epsilon}+\tau^{\epsilon}\nabla^T u ^{\epsilon})\|_{L^2}\notag\\
&+\|\nabla^{\ell}\big(\eta^{\epsilon}(\nabla  u ^{\epsilon}+\nabla^T u ^{\epsilon})\big)\|_{L^2}
+\|\Delta\nabla^{\ell} \tau^{\epsilon}\|_{L^2}+\|\nabla^{\ell}\tau^{\epsilon}\|_{L^2}\label{4-s29}\\
\lesssim& \|\tau ^{\epsilon}\|_{H^2}\| u ^{\epsilon}\|_{H^2}+(\|\nabla\eta ^{\epsilon}\|_{H^1}+\| \eta ^{\epsilon}\|_{L^\infty})\|\nabla u ^{\epsilon}\|_{H^1}+\|\tau ^{\epsilon}\|_{H^3}\notag\\
\lesssim& \delta.\notag
\end{align}
Combining \eqref{4-s28} with \eqref{4-s29}, we can get the inequality \eqref{4-4}.

Secondly, it is worth noticing that we cannot directly estimate $\partial_t\phi^{\epsilon}$ and $\partial_t  u ^{\epsilon}$ by following the above method due to the singularity in the first two equations of \eqref{ma1}. Therefore, the main difficulty in subsequent estimation is to cancel out the singularity. First of all, applying derivatives $\partial_t$ to the first two equations of \eqref{ma1}, we get
\beq\label{4-ma3}
\left\{
\begin{split}
&\partial_{tt}\phi^{\epsilon}
+ u ^{\epsilon}\cdot\nabla\partial_t\phi^{\epsilon}
+\partial_t u ^{\epsilon}\cdot\nabla\phi^{\epsilon}
+\partial_t\phi^{\epsilon}\mathrm{div} u ^{\epsilon}
+\phi^{\epsilon}\mathrm{div}\partial_t  u ^{\epsilon}
+\frac{1}{\epsilon}\mathrm{div}\partial_t u ^{\epsilon}=0,\\
&\partial_{tt}  u ^{\epsilon}
+ \partial_t u ^{\epsilon}\cdot \nabla  u ^{\epsilon}
+  u ^{\epsilon}\cdot \nabla \partial_t u ^{\epsilon}
+\frac{1}{\epsilon}\partial_t\big(\frac{P'(\rho^{\epsilon})}{\rho^{\epsilon}}\big) \nabla\phi^{\epsilon}
+\frac{1}{\epsilon}\frac{P'(\rho^{\epsilon})}{\rho^{\epsilon}} \nabla\partial_t\phi^{\epsilon}=\partial_t\big(\frac{1}{\rho^{\epsilon}}\big)B+\frac{1}{\rho^{\epsilon}}\partial_tB,
\end{split}
\right.
\eeq
where~$B=\mu_1\Delta  u ^{\epsilon}+\mu_2\nabla \mathrm{div} u ^{\epsilon}+\frac{\beta}{k}\mathrm{div}\tau^{\epsilon}
-\nabla\big(\beta(L-1)\eta^{\epsilon}+\bar{\mathfrak{z}}(\eta^{\epsilon})^{2}\big)$. Then, for $\ell=0,1$, applying the derivative operator $\nabla^\ell$ to the equation \eqref{4-ma3}$_1$, and taking the inner
product with $P'(\rho^{\epsilon})\nabla^{\ell}\partial_t\phi^{\epsilon}$, we find that that for some small constant $\theta_3>0$ independent of $\epsilon$,
\begin{align}
&\frac{1}{2}\frac{d}{dt}\|\partial_t\phi^{\epsilon}\|_{H^1_{P'(\rho^{\epsilon})}}^2
+\frac{1}{\epsilon}\sum\limits_{\ell=0}^1\langle\mathrm{div}\nabla^{\ell}\partial_t u ^{\epsilon},
P'(\rho^{\epsilon})\nabla^{\ell}\partial_t\phi^{\epsilon}\rangle\label{4-s32}\\
=&\frac{1}{2}\sum\limits_{\ell=0}^1\langle\partial_t P'(\rho^{\epsilon}),|\nabla^{\ell}\partial_t\phi^{\epsilon}|^2\rangle
-\sum\limits_{\ell=0}^1\langle \nabla^{\ell}( u ^{\epsilon}\cdot\nabla\partial_t\phi^{\epsilon})
,P'(\rho^{\epsilon})\nabla^{\ell}\partial_t\phi^{\epsilon}\rangle\notag\\
&-\sum\limits_{\ell=0}^1\langle \nabla^{\ell}(\partial_t u ^{\epsilon}\cdot\nabla\phi^{\epsilon})
,P'(\rho^{\epsilon})\nabla^{\ell}\partial_t\phi^{\epsilon}\rangle
-\sum\limits_{\ell=0}^1\langle \nabla^{\ell}(\mathrm{div}  u ^{\epsilon}\partial_t\phi^{\epsilon})
,P'(\rho^{\epsilon})\nabla^{\ell}\partial_t\phi^{\epsilon}\rangle\notag\\
&
-\sum\limits_{\ell=0}^1\langle \nabla^{\ell}(\mathrm{div}\partial_t  u ^{\epsilon}\phi^{\epsilon})
,P'(\rho^{\epsilon})\nabla^{\ell}\partial_t\phi^{\epsilon}\rangle\notag\\
\leq& C_{\theta_3}\big(\|\partial_t\phi^{\epsilon}\|_{H^1_{P'(\rho^{\epsilon})}}^2
+\|\partial_t  u ^{\epsilon}\|_{H^1_{\rho^{\epsilon}}}^2\big)
+\theta_3\|\nabla\mathrm{div}\partial_t  u ^{\epsilon}\|_{L^2}^2,\notag
\end{align}
where the last inequality is derived based on \eqref{clc-02}, \eqref{clc-03}, \eqref{4-s14}, integration by parts, H\"older inequality, Sobolev inequality and Cauchy inequality. Similarly, for $\ell=0,1$, applying the derivative operator $\nabla^\ell$ to the equation \eqref{4-ma3}$_2$, and taking the inner product with $\rho^{\epsilon}\nabla^{\ell}\partial_t  u ^{\epsilon}$, we can obtain
\begin{align}
&\frac{1}{2}\frac{d}{dt}\|\partial_t  u ^{\epsilon}\|_{H^1_{\rho^{\epsilon}}}^2
+\mu_1\|\nabla\partial_t u ^{\epsilon}\|_{H^1}^2
+\mu_2\|\mathrm{div}\partial_t u ^{\epsilon}\|_{H^1}^2
+\frac{1}{\epsilon}\sum\limits_{\ell=0}^1
\langle\nabla^{\ell}\partial_t u ^{\epsilon},
P'(\rho^{\epsilon})\nabla^{\ell+1}\partial_t\phi^{\epsilon}\rangle\notag\\
=&-\sum\limits_{\ell=0}^1\langle [\nabla^{\ell}, u ^{\epsilon}\cdot\nabla]\partial_t  u ^{\epsilon},\rho^{\epsilon}\nabla^{\ell}\partial_t  u ^{\epsilon}\rangle
-\frac{1}{\epsilon}\sum\limits_{\ell=0}^1
\langle[\nabla^{\ell},\frac{P'(\rho^{\epsilon})}{\rho^{\epsilon}}]\nabla
\partial_t\phi^{\epsilon},\rho^{\epsilon}\nabla^{\ell}\partial_t u ^{\epsilon}
\rangle\label{4-s25}\\
&-\sum\limits_{\ell=0}^1\langle \nabla^{\ell}(\partial_t u ^{\epsilon}\cdot\nabla  u ^{\epsilon}),\rho^{\epsilon}\nabla^{\ell}\partial_t  u ^{\epsilon}\rangle
-\frac{1}{\epsilon}\sum\limits_{\ell=0}^1\langle \nabla^{\ell}(\partial_t\frac{P'(\rho^{\epsilon})}{\rho^{\epsilon}}
\nabla\phi^{\epsilon}),\rho^{\epsilon}\nabla^{\ell}\partial_t  u ^{\epsilon}\rangle\notag\\
&+\sum\limits_{\ell=0}^1\langle \nabla^{\ell}\big(\partial_t( \frac{1}{\rho^{\epsilon}})
B\big),\rho^{\epsilon}\nabla^{\ell}\partial_t  u ^{\epsilon}\rangle
+\sum\limits_{\ell=0}^1\langle[ \nabla^{\ell},\frac{1}{\rho^{\epsilon}}
]\partial_tB,\rho^{\epsilon}\nabla^{\ell}\partial_t  u ^{\epsilon}\rangle\notag\\
&+\sum\limits_{\ell=0}^1\langle \nabla^{\ell}\partial_t\left(\frac{\beta}{k}\mathrm{div}\tau^{\epsilon}
-\nabla\big(\beta(L-1)\eta^{\epsilon}+\bar{\mathfrak{z}}(\eta^{\epsilon})^{2}\big)\right),
\nabla^{\ell}\partial_t  u ^{\epsilon}\rangle,\notag
\end{align}
where we have used the fact that
\[\frac{1}{2}\langle\partial_t \rho^{\epsilon},|\nabla^{\ell}\partial_t u ^{\epsilon}|^2\rangle
=\langle  u ^{\epsilon}\cdot\nabla^{\ell+1}\partial_t  u ^{\epsilon},\rho^{\epsilon}\nabla^{\ell}\partial_t  u ^{\epsilon}\rangle.\]
Next, we are going to deal with the terms on the right-hand side of \eqref{4-s25}. For the first line on the right-hand
side of \eqref{4-s25}, noticing that it is equal to $0$ when $\ell=0$, we can use H\"older inequality, Sobolev inequality, Cauchy inequality, \eqref{clc-02} and \eqref{clc-03} to obtain
\begin{align}
&-\sum\limits_{\ell=0}^1\langle [\nabla^{\ell}, u ^{\epsilon}\cdot\nabla]\partial_t  u ^{\epsilon},\rho^{\epsilon}\nabla^{\ell}\partial_t  u ^{\epsilon}\rangle
-\frac{1}{\epsilon}\sum\limits_{\ell=0}^1
\langle[\nabla^{\ell},(\frac{P'(\rho^{\epsilon})}{\rho^{\epsilon}}-a\gamma)]\nabla
\partial_t\phi^{\epsilon},\rho^{\epsilon}\nabla^{\ell}\partial_t u ^{\epsilon}
\rangle\notag\\
\lesssim&\|\nabla  u ^{\epsilon}\|_{L^\infty}\|\nabla\partial_t  u ^{\epsilon}\|_{L^2}\|\rho^{\epsilon}\|_{L^\infty}
\|\nabla\partial_t  u ^{\epsilon}\|_{L^2}
+\|\rho^{\epsilon}\|_{L^\infty}\|\nabla\phi^{\epsilon}\|_{L^\infty}\|\nabla\partial_t \phi^{\epsilon}\|_{L^2}
\|\nabla\partial_t  u ^{\epsilon}\|_{L^2}\label{4-s26}\\
\lesssim&\|\nabla\partial_t \phi^{\epsilon}\|_{L^2_{P'(\rho^{\epsilon})}}^2
+\|\nabla\partial_t  u ^{\epsilon}\|_{L^2_{\rho^{\epsilon}}}^2\notag.
\end{align}
Similarly, the second line on the right-hand side of \eqref{4-s25} can be controlled as follows
\begin{align}
&-\sum\limits_{\ell=0}^1\langle \nabla^{\ell}(\partial_t u ^{\epsilon}\cdot\nabla  u ^{\epsilon}),\rho^{\epsilon}\nabla^{\ell}\partial_t  u ^{\epsilon}\rangle
-\frac{1}{\epsilon}\sum\limits_{\ell=0}^1\langle \nabla^{\ell}(\partial_t\frac{P'(\rho^{\epsilon})}{\rho^{\epsilon}}
\nabla\phi^{\epsilon}),\rho^{\epsilon}\nabla^{\ell}\partial_t  u ^{\epsilon}\rangle\label{4-s27}\\
\lesssim&\|\nabla  u ^{\epsilon}\|_{L^\infty}\|\partial_t  u ^{\epsilon}\|_{H^1}\|\rho^{\epsilon}\|_{L^\infty}
\|\partial_t  u ^{\epsilon}\|_{H^1}
+\|\nabla^2  u ^{\epsilon}\|_{L^3}\|\partial_t  u ^{\epsilon}\|_{L^6}\|\rho^{\epsilon}\|_{L^\infty}
\|\nabla\partial_t  u ^{\epsilon}\|_{L^2}\notag\\
&+\|\nabla \phi^{\epsilon}\|_{L^\infty}\|\partial_t \phi^{\epsilon}\|_{H^1}\|\rho^{\epsilon}\|_{L^\infty}
\|\partial_t  u ^{\epsilon}\|_{H^1}
+\|\nabla^2 \phi^{\epsilon}\|_{L^3}\|\partial_t \phi^{\epsilon}\|_{L^6}\|\rho^{\epsilon}\|_{L^\infty}
\|\nabla\partial_t  u ^{\epsilon}\|_{L^2}\notag\\
\lesssim&\|\partial_t \phi^{\epsilon}\|_{H^1_{P'(\rho^{\epsilon})}}^2
+\|\partial_t  u ^{\epsilon}\|_{H^1_{\rho^{\epsilon}}}^2\notag.
\end{align}
By the virtue of ~\eqref{clc-02}, \eqref{clc-03} and \eqref{4-4}, the third line on the right-hand side of \eqref{4-s25} can be estimated as
\begin{align}
&\sum\limits_{\ell=0}^1\langle \nabla^{\ell}\big(\partial_t(\frac{1}{\rho^{\epsilon}})
B\big),\rho^{\epsilon}\nabla^{\ell}\partial_t  u ^{\epsilon}\rangle
+\sum\limits_{\ell=0}^1\langle[ \nabla^{\ell},(\frac{1}{\rho^{\epsilon}}-1)
]\partial_tB,\rho^{\epsilon}\nabla^{\ell}\partial_t  u ^{\epsilon}\rangle\notag\\
\lesssim&\|\frac{1}{\rho^{\epsilon}}\|_{L^\infty}^2\|\partial_t \phi^{\epsilon}\|_{L^2}\|B\|_{L^\infty}\|\rho^{\epsilon}\|_{L^\infty}
\|\partial_t  u ^{\epsilon}\|_{L^2}\label{4-s30}\\
&+\|B\|_{L^\infty}\big(\|\frac{1}{\rho^{\epsilon}}\|_{L^\infty}^2\|\nabla\partial_t \phi^{\epsilon}\|_{L^2}+\|\frac{1}{\rho^{\epsilon}}\|_{L^\infty}^3\|\nabla \phi^{\epsilon}\|_{L^\infty}\|\partial_t \phi^{\epsilon}\|_{L^2}\big)\|\rho^{\epsilon}\|_{L^\infty}
\|\nabla\partial_t  u ^{\epsilon}\|_{L^2}\notag\\
&
+\|\frac{1}{\rho^{\epsilon}}\|_{L^\infty}^2\big(\|\nabla B\|_{L^2}\|\partial_t \phi^{\epsilon}\|_{L^3}\|\rho^{\epsilon}\|_{L^\infty}
\|\nabla\partial_t  u ^{\epsilon}\|_{L^6}+\|\partial_t B\|_{L^2}\|\nabla \phi^{\epsilon}\|_{L^\infty}\|\rho^{\epsilon}\|_{L^\infty}
\|\nabla\partial_t  u ^{\epsilon}\|_{L^2}\big)\notag\\
\leq&C_{\theta_3}(\|\partial_t \phi^{\epsilon}\|_{H^1_{P'(\rho^{\epsilon})}}^2
+\|\partial_t  u ^{\epsilon}\|_{H^1_{\rho^{\epsilon}}}^2+1)
+\theta_3\|\nabla^2\partial_t  u ^{\epsilon}\|_{L^2}^2\notag,
\end{align}
where we have used H\"older inequality, Sobolev inequality, Cauchy inequality and the following inequality
\begin{align*}
\|\partial_tB\|_{L^2}\lesssim\|\nabla^2 \partial_t  u ^{\epsilon}\|_{L^2}
+\|\nabla\partial_t\eta^{\epsilon}\|_{L^2}(\|\nabla\eta^{\epsilon}\|_{H^1}+1)
+\|\nabla\partial_t\tau^{\epsilon}\|_{L^2}
\lesssim \|\nabla^2 \partial_t  u ^{\epsilon}\|_{L^2}+1.
\end{align*}
For the last line on the right-hand side of \eqref{4-s25}, thanks to integration by parts, it is not hard to get
\begin{align}
&\sum\limits_{\ell=0}^1\langle \nabla^{\ell}\partial_t\left(\frac{\beta}{k}\mathrm{div}\tau^{\epsilon}
-\nabla\big(\beta(L-1)\eta^{\epsilon}+\bar{\mathfrak{z}}(\eta^{\epsilon})^{2}\big)\right),
\nabla^{\ell}\partial_t  u ^{\epsilon}\rangle\label{4-s31}\\
=&\langle\partial_t\left(\frac{\beta}{k}\mathrm{div}\tau^{\epsilon}
-\nabla\big(\beta(L-1)\eta^{\epsilon}+\bar{\mathfrak{z}}(\eta^{\epsilon})^{2}\big)\right),
\partial_t  u ^{\epsilon}-\Delta\partial_t  u ^{\epsilon}\rangle\notag\\
\lesssim &\big(\|\nabla\partial_t\tau^{\epsilon}\|_{L^2}
+\|\nabla\partial_t\eta^{\epsilon}\|_{L^2}
+\|\nabla\eta^{\epsilon}\|_{H^1}
\|\nabla\partial_t\eta^{\epsilon}\|_{L^2}\big)(\|\partial_t  u ^{\epsilon}\|_{L^2(\rho^{\epsilon})}+\|\nabla^2\partial_t  u ^{\epsilon}\|_{L^2})\notag\\
\leq &C_{\theta_3}(\|\partial_t  u ^{\epsilon}\|_{L^2_{\rho^{\epsilon}}}^2+1)
+\theta_3\|\nabla^2\partial_t  u ^{\epsilon}\|_{L^2}^2.\notag
\end{align}
Then, combining \eqref{4-s32} with \eqref{4-s25}-\eqref{4-s31}, choosing $\theta_3$ small enough, we find that there is a positive constant $C_3$ independent
of $\epsilon$, such that
\begin{align}
&\frac{d}{dt}(\|\partial_t \phi^{\epsilon}\|_{H^1_{P'(\rho^{\epsilon})}}^2+
\|\partial_t  u ^{\epsilon}\|_{H^1_{\rho^{\epsilon}}}^2)
+\mu_1\|\nabla\partial_t u ^{\epsilon}\|_{H^1}^2
+\mu_2\|\mathrm{div}\partial_t u ^{\epsilon}\|_{H^1}^2\notag\\
\lesssim& 1+\|\partial_t \phi^{\epsilon}\|_{H^1_{P'(\rho^{\epsilon})}}^2+
\|\partial_t  u ^{\epsilon}\|_{H^1_{\rho^{\epsilon}}}^2
+\frac{1}{\epsilon}\sum\limits_{\ell=0}^1\langle\nabla^{\ell}\partial_t u ^{\epsilon},
\nabla P'(\rho^{\epsilon})\nabla^{\ell}\partial_t\phi^{\epsilon}\rangle\label{4-s33}\\
\lesssim& 1+\|\partial_t \phi^{\epsilon}\|_{H^1_{P'(\rho^{\epsilon})}}^2+
\|\partial_t  u ^{\epsilon}\|_{H^1_{\rho^{\epsilon}}}^2
+\|P''(\rho^{\epsilon})\|_{L^\infty}\|\nabla\phi^{\epsilon}\|_{L^\infty}
\|\partial_t  u ^{\epsilon}\|_{H^1_{\rho^{\epsilon}}}\|\partial_t \phi^{\epsilon}\|_{H^1_{P'(\rho^{\epsilon})}}\notag\\
\leq&C_3(1+\|\partial_t \phi^{\epsilon}\|_{H^1_{P'(\rho^{\epsilon})}}^2+
\|\partial_t  u ^{\epsilon}\|_{H^1_{\rho^{\epsilon}}}^2)\notag,
\end{align}
where we have used the following relation to deal with the singular terms that appears on the left side of \eqref{4-s32} and \eqref{4-s25}.
\begin{align*}
&\frac{1}{\epsilon}\sum\limits_{\ell=0}^1\langle\mathrm{div}\nabla^{\ell}\partial_t u ^{\epsilon},
P'(\rho^{\epsilon})\nabla^{\ell}\partial_t\phi^{\epsilon}\rangle
+\frac{1}{\epsilon}\sum\limits_{\ell=0}^1
\langle\nabla^{\ell}\partial_t u ^{\epsilon},
P'(\rho^{\epsilon})\nabla^{\ell+1}\partial_t\phi^{\epsilon}\rangle\notag\\
=&-\frac{1}{\epsilon}\sum\limits_{\ell=0}^1\langle\nabla^{\ell}\partial_t u ^{\epsilon},
\nabla P'(\rho^{\epsilon})\nabla^{\ell}\partial_t\phi^{\epsilon}\rangle\notag.
\end{align*}
Further, by virtue of Gronwall's inequality, \eqref{4-s33} implies that for any fixed $T>0$,
\begin{align}
\|\partial_t \phi^{\epsilon}\|_{H^1_{P'(\rho^{\epsilon})}}^2+
\|\partial_t  u ^{\epsilon}\|_{H^1_{\rho^{\epsilon}}}^2
\leq \big(1+\|\partial_t \phi^{\epsilon}(0)\|_{H^1_{P'(\rho^{\epsilon}_0)}}^2+
\|\partial_t  u ^{\epsilon}(0)\|_{H^1_{\rho^{\epsilon}_0}}^2\big)\exp(C_3T)\label{4-s34},
\end{align}
holds for all $t\in[0,T]$ and $\epsilon$. It can be seen from inequality \eqref{4-s34} that in order to estimate $\|\partial_t \phi^{\epsilon}\|_{H^1_{P'(\rho^{\epsilon}_0)}}$ and $
\|\partial_t  u ^{\epsilon}\|_{H^1_{\rho^{\epsilon}_0}}$, we have to first obtain the estimates of $\|\partial_t \phi^{\epsilon}(0)\|_{H^1_{P'(\rho^{\epsilon}_0)}}$ and $
\|\partial_t  u ^{\epsilon}(0)\|_{H^1_{\rho^{\epsilon}_0}}$. In fact, recalling the constraints of the initial data \eqref{initial2} and \eqref{clc-05}, we get
\begin{align}\label{4-s75}
\|\partial_t\phi^{\epsilon}(0)\|_{H^1}\lesssim \| \nabla  u ^{\epsilon}_0\|_{H^1}\| \nabla\phi^{\epsilon}_0\|_{H^1}
+\frac{1}{\epsilon}\|\mathrm{div} u ^{\epsilon}_0\|_{H^1}
\lesssim \delta+1,
\end{align}
and
\begin{align}
\|\partial_t  u ^{\epsilon}(0)\|_{H^1}\lesssim& \| \nabla  u ^{\epsilon}_0\|_{H^1}^2
+\frac{1}{\epsilon}\|\frac{P'(\rho^{\epsilon}_0)}{\rho^{\epsilon}_0} \nabla\phi^{\epsilon}_0 \|_{H^1}
+(1+\|\nabla\phi^{\epsilon}_0\|_{H^1}^2)
\|\nabla\eta^{\epsilon}_0\|_{H^1}\notag\\
&+(1+\|\nabla\phi^{\epsilon}_0\|_{H^1}^2)
\big(\|\nabla\eta^{\epsilon}_0\|_{H^1}^2
+\|\nabla^2 u ^{\epsilon}_0\|_{H^1}+\|\nabla\tau^{\epsilon}_0\|_{H^1}\big)\label{4-s76}\\[1mm]
\lesssim& \,\delta+1,\notag
\end{align}
where we have used the fact that
\begin{align*}
\frac{1}{\epsilon}\|\frac{P'(\rho^{\epsilon}_0)}{\rho^{\epsilon}_0} \nabla\phi^{\epsilon}_0 \|_{H^1}
\lesssim&\frac{1}{\epsilon}\bigg(\|\frac{P'(\rho^{\epsilon}_0)}{\rho^{\epsilon}_0}\|_{L^\infty}
\|\nabla\phi^{\epsilon}_0 \|_{L^2}+\|\nabla\frac{P'(\rho^{\epsilon}_0)}{\rho^{\epsilon}_0}\|_{L^3}
\|\nabla\phi^{\epsilon}_0 \|_{L^6}
+\|\frac{P'(\rho^{\epsilon}_0)}{\rho^{\epsilon}_0}\|_{L^\infty}
\|\nabla^2\phi^{\epsilon}_0 \|_{L^2}\bigg)\\
\lesssim&\frac{1}{\epsilon}\|\nabla\phi^{\epsilon}_0 \|_{H^1}
\leq C.
\end{align*}
Therefore, \eqref{4-s34} together with \eqref{4-s75} and \eqref{4-s76} yields
\begin{align*}
\|\partial_t\phi^{\epsilon}\|_{L^{\infty}(0,T;H^1_{P'(\rho^{\epsilon})})}+
\|\partial_t  u ^{\epsilon}\|_{L^{\infty}(0,T;H^1_{\rho^{\epsilon}})}
\leq C(T).
\end{align*}

Finally, we conclude from the first two equations of \eqref{ma1}, the uniform bounds \eqref{clc-02}, \eqref{clc-03} and \eqref{4-5} that
\begin{align*}
\frac{1}{\epsilon}\|\mathrm{div} u ^{\epsilon}\|_{H^1}
\lesssim \|\partial_t\phi^{\epsilon}\|_{H^1}+\| u ^{\epsilon}\cdot\nabla\phi^{\epsilon}\|_{H^1}
+\|\phi^{\epsilon}\mathrm{div} u ^{\epsilon}\|_{H^1}\leq C(T)
\end{align*}
and
\begin{align*}
\frac{1}{\epsilon}\|P'(\rho^{\epsilon})\nabla \phi^{\epsilon}\|_{H^1}
\lesssim & \|\rho^{\epsilon}\partial_t u ^{\epsilon}\|_{H^1}+\|\rho^{\epsilon} u ^{\epsilon}\cdot \nabla u ^{\epsilon}\|_{H^1}
+\|\nabla\eta^{\epsilon}\|_{H^1}\\
&+\|\nabla\eta^{\epsilon}\|_{H^1}^{2}
+\|\nabla^2  u ^{\epsilon}\|_{H^1}+\|\nabla \tau^{\epsilon}\|_{H^1}\\
\leq& C(T).
\end{align*}
Combining the last two inequalities, we obtain \eqref{4-6} and complete the proof of this proposition.
\endpf

\section{Proof of Theorem \ref{Theorem2}}\label{sec4}
Based on \eqref{clc-02}, there exists a subsequence of $(\phi^\epsilon,u^\epsilon,\eta^\epsilon,\tau^\epsilon)$, which will be still denoted as $(\phi^\epsilon,u^\epsilon,\eta^\epsilon,\tau^\epsilon)$ for convenience, with a limit $(\phi,u,\eta,\tau)$, such that for any fixed $T>0$,
\begin{equation}\label{s39}
\begin{split}
  &(\phi^\epsilon,u^\epsilon,\eta^\epsilon-1,\tau^\epsilon)\rightarrow(\phi,u,\eta-1,\tau) \;\;\text{weakly-$*$ in}\;\; L^{\infty}(\mathbb{R}^+;H^3),\\
  &(\nabla u^\epsilon,\nabla\eta^\epsilon,\tau^\epsilon, \nabla\tau^\epsilon)\rightarrow(\nabla u,\nabla\eta,\tau, \nabla\tau) \;\;\text{weakly in}\;\; L^2(\mathbb{R}^+;H^3),
 \end{split}
\end{equation}
as $\epsilon\rightarrow0$. Then, owing to \eqref{clc-02}, \eqref{clc-04}, \eqref{clc-06} and Lemma \ref{aubin}, we deduce that
\begin{align}
(\phi^\epsilon,u^\epsilon,\eta^\epsilon-1,\tau^\epsilon)
\rightarrow(\phi,u,\eta-1,\tau)\;\; \text{strongly in}\;\;\mathcal{C}([0,T]; H^2_{\text{loc}}),\label{s35}
\end{align} for any fixed $T>0$, as $\epsilon\rightarrow0$.
In what follows, we show that $(u,\eta,\tau)$ is the unique strong solution to the incompressible Oldroyd-B model \eqref{limit1} by virtue of the corresponding compressible model \eqref{ma1} as well as the limits \eqref{s39} and \eqref{s35}.

Firstly, we prove that $(u,\eta,\tau)$ satisfies the equations of $u$ stated in \eqref{limit1}. To show this, we rewritten the system \eqref{ma1}$_2$ as follows
\begin{align}
\partial_t u^{\epsilon}+ u^{\epsilon}\cdot \nabla u^{\epsilon}
+\epsilon\phi^{\epsilon}\partial_t u^{\epsilon}+ \epsilon\phi^{\epsilon}u^{\epsilon}\cdot \nabla u^{\epsilon}
+&\frac{1}{\epsilon}P'(\rho^{\epsilon}) \nabla\phi^{\epsilon} +\nabla\big(\beta(L-1)\eta^{\epsilon}+\bar{\mathfrak{z}}(\eta^{\epsilon})^{2}\big)\label{s1-}\\
&=\mu_1\Delta u^{\epsilon}+\mu_2\nabla \mathrm{div}u^{\epsilon}+\frac{\beta}{k}\mathrm{div} \tau^{\epsilon}.\notag
\end{align}
Then, we  analyze the convergence of each term on the equation \eqref{s1-}. Thanks to \eqref{clc-03}, \eqref{clc-06}, \eqref{clc-07} and \eqref{s35} again, there exists a function $v\in L^{\infty}(0,T;H^1)$, such that
\begin{equation}\label{s38-}
\begin{split}
&\partial_t u^{\epsilon}\rightarrow \partial_t u\;\;\text{weakly-$*$ in}\;\; L^\infty(0,T;H^1),\\
&\frac{1}{\epsilon}P'(\rho^\epsilon)\nabla\phi^\epsilon
\rightarrow v\;\;\text{weakly-$*$ in}\;\; L^{\infty}(0,T;H^1),\;\;\text{for any fixed}\;\;T>0,
\end{split}
\end{equation}
as $\epsilon\rightarrow0$. In addition, we use the uniform bounds \eqref{clc-02} and \eqref{clc-06} to give
\begin{align*}
&\|\phi^{\epsilon}\partial_t u^{\epsilon}\|_{L^\infty(0,T;H^1) }
\leq C\|\phi^{\epsilon}\|_{L^\infty(0,T;H^3)}
\|\partial_t u^{\epsilon}\|_{L^\infty(0,T;H^1)}
\leq C(T),\\
&\|\phi^{\epsilon}u^{\epsilon}\cdot\nabla u ^{\epsilon}\|_{L^\infty(\mathbb{R}^+;H^2)}
\leq C\|\phi^{\epsilon}\|_{L^\infty(\mathbb{R}^+;H^3)}
\|u^{\epsilon}\|_{L^\infty(\mathbb{R}^+;H^3)}^2
\leq C,
\end{align*}
which implies
\begin{equation}\label{s38}
\begin{split}
&\epsilon \phi^{\epsilon}\partial_t u^{\epsilon}\rightarrow 0\;\;\text{strongly in}\;\; L^\infty(0,T;H^1),\\
&\epsilon \phi^{\epsilon}u^{\epsilon}\cdot\nabla u^{\epsilon}\rightarrow 0\;\;\text{strongly in}\;\; L^\infty(\mathbb{R}^+;H^2), \;\;\text{for any fixed}\;\;T>0,
\end{split}
\end{equation}
as $\epsilon\rightarrow0$. Moreover, by virtue of the convergence \eqref{s35}, it is easy to verify that
\begin{equation}\label{s41}
\begin{split}
&\nabla\big(\beta(L-1)\eta^{\epsilon}+\bar{\mathfrak{z}}(\eta^{\epsilon})^{2}\big)\rightarrow \nabla\big(\beta(L-1)\eta+\bar{\mathfrak{z}}(\eta)^{2}\big)\;\;\text{strongly in}\;\; \mathcal{C}(0,T;H^1_{\text{loc}}),\\
&\mu_1\Delta u^{\epsilon}+\mu_2\nabla \mathrm{div}u^{\epsilon}\rightarrow \mu_1\Delta u+\mu_2\nabla \mathrm{div}u\;\;\text{strongly in}\;\; \mathcal{C}(0,T;L^2_{\text{loc}}),\\
&u^{\epsilon}\cdot \nabla u^{\epsilon}\rightarrow u\cdot \nabla u\;\;\text{strongly in}\;\; \mathcal{C}(0,T;H^1_{\text{loc}}),\\
&\mathrm{div}\tau^{\epsilon}\rightarrow\mathrm{div}\tau \;\;\text{strongly in}\;\; \mathcal{C}(0,T;H^1_{\text{loc}}), \;\;\text{for any fixed}\;\;T>0,
\end{split}
\end{equation}
as $\epsilon\rightarrow0$. Let $\varphi_1(x,t)$ be a smooth test function of \eqref{s1-} with compact supports in $\mathbb{R}^3\times(0,T)$ and satisfying the divergence free condition $\mathrm{div}\varphi_1=0$, then
\begin{align*}
&\int_0^{T}\int_{\mathbb{R}^3}\big(\partial_t u^{\epsilon}+ u^{\epsilon}\cdot \nabla u^{\epsilon}
+\epsilon\phi^{\epsilon}\partial_t u^{\epsilon}+ \epsilon\phi^{\epsilon}u^{\epsilon}\cdot \nabla u^{\epsilon}\big)\cdot\varphi_1\;\mathrm{d}x\mathrm{d}t\\
&+\int_0^{T}\int_{\mathbb{R}^3}\bigg(\nabla\big(\beta(L-1)\eta^{\epsilon}+\bar{\mathfrak{z}}(\eta^{\epsilon})^{2}\big)
-\mu_1\Delta u^{\epsilon}-\mu_2\nabla \mathrm{div}u^{\epsilon}-\frac{\beta}{k}\mathrm{div} \tau^{\epsilon}\bigg)\cdot\varphi_1\;\mathrm{d}x\mathrm{d}t\\
=&-\frac{1}{\epsilon^2}\int_0^{T}\int_{\mathbb{R}^3}\varphi_1\cdot \nabla P(\rho^{\epsilon})\;\mathrm{d}x\mathrm{d}t=0.
\end{align*}
Sending $\epsilon\rightarrow0$ in the above equation, and using \eqref{s38-}-\eqref{s41}, we get
\begin{align*}
\partial_t u+ u\cdot \nabla u +\nabla\big(\beta(L-1)\eta+\bar{\mathfrak{z}}(\eta)^{2}\big)
-\mu_1\Delta u-\mu_2\nabla \mathrm{div}u-\frac{\beta}{k}\mathrm{div} \tau=-\nabla
\pi_1,
\end{align*}
for some $\pi_1\in L^{\infty}(0,T;H^2)$. Observing the limit \eqref{s38-}, we can further deduce that $v=\nabla \pi_1$ in $\mathbb{R}^3\times [0,T]$, which implies
\begin{align*}
\frac{1}{\epsilon}P'(\rho^\epsilon)\nabla\phi^\epsilon
\rightarrow \nabla \pi_1\;\;\text{weakly-$*$ in}\;\; L^{\infty}(0,T;H^1),\;\;\text{for any fixed}\;\;T>0,
\end{align*}as $\epsilon\rightarrow0$.
On the other hand, the initial condition \eqref{clc-08} and the convergence \eqref{s39} yield that $u(x,0)=u_0(x)$ for a.e. $x\in \mathbb{R}^3$. Furthermore, from the uniform bound \eqref{clc-07}, it is easily to show that $\mathrm{div}u^\epsilon\rightarrow 0$ strongly in $L^\infty(0,T;H^1)$ for any fixed $T>0$, as $\epsilon\rightarrow0$, which, combining with the convergence \eqref{s35}, implies that
\begin{align*}
\mathrm{div}u=0\;\;\text{for a.e.}\;\; (x,t)\in \mathbb{R}^3\times [0,T].
\end{align*}
In summery, we have shown that $u(x,t)\in L^\infty(\mathbb{R}^+;H^3)$ and $\nabla u(x,t)\in L^2(\mathbb{R}^+;H^3)$ satisfy
\begin{align*}
\partial_t u+ u\cdot \nabla u-\mu_1\Delta u+\nabla\pi
&=\frac{\beta}{k}\mathrm{div} \tau,\\
\mathrm{div}u&=0,
\end{align*}
with initial data
\[u|_{t=0}=u_0(x),\]
where $\pi=\pi_1+\beta(L-1)\eta+\bar{\mathfrak{z}}(\eta)^{2}$.

Similarly, we can also deduce $(\eta-1,\tau)(x,t)\in L^\infty(\mathbb{R}^+;H^3)$ and $(\nabla\eta,\tau,\nabla\tau) (x,t)\in L^2(\mathbb{R}^+;H^3)$ satisfy the following equation
\begin{align*}
&\partial_t\eta+\mathrm{div} (\eta  u)=\nu\Delta \eta,\\
&\partial_t\tau+\mathrm{div} (u\tau)-(\nabla u
\tau+\tau\nabla^Tu)-k\eta(\nabla u+\nabla^Tu)
=\nu\Delta \tau-\frac{A_0}{2}\tau,
\end{align*}
with initial data
\[(\eta, \tau)|_{t=0}=(\eta_0,\tau_0)(x).\]

Therefore, the limit functions $(u,\eta,\tau)(x,t)$ satisfy the system \eqref{limit1} with the initial data \eqref{limit-initial}. In addition, the uniform bound \eqref{clc-09} can be deduced from the uniform bound \eqref{clc-02} and the lower semi-continuity of norms. Moreover, the regularity of the solution is good enough to ensure the uniqueness.

In a word, we have proved that $(\phi^\epsilon,u^\epsilon,\eta^\epsilon,\tau^\epsilon)$ converges to $(\phi,u,\eta,\tau)$ as $\epsilon\rightarrow0$ and $(u,\eta,\tau)$ is the unique global strong solution to the incompressible Oldroyd-B model \eqref{limit1}. In order to finish the proof of Theorem \ref{Theorem2}, what left is to show the convergence rate by using the relative entropy method. To achieve this, we need the following lemma.

\begin{lemma}\label{lemma1}(\cite{GJLLT}) Under the same assumptions in Theorem \ref{Theorem2}, we have
\begin{align*}
\|\sqrt{\rho^\epsilon}-1\|_{L^2}\lesssim \epsilon\langle\Pi^\epsilon,1\rangle^\frac{1}{2},\;\;\;
\|\sqrt{\rho^\epsilon_0}-1\|_{L^2}\lesssim \epsilon\langle\Pi^\epsilon_0,1\rangle^\frac{1}{2}
\lesssim \epsilon^{1+\frac{\alpha_0}{2}},
\end{align*}
where $\Pi^{\epsilon}=\frac{1}{\epsilon^2}\frac{a}{\gamma-1}[(\rho^{\epsilon})^{\gamma}-\gamma(\rho^{\epsilon}-1)-1]$ and $\alpha_0>0$ is a constant independent of $\epsilon$.
\end{lemma}

In the sequel, we begin the proof of the convergence rate stated in Theorem \ref{Theorem2}. Firstly, multiplying the three equations \eqref{ma1}$_2$, \eqref{ma1}$_3$ and \eqref{ma1}$_4$ by $\rho^\epsilon u^\epsilon$, $\eta^\epsilon-1$ and  $\tau^\epsilon$, respectively, integrating the resulting equations over $\mathbb{R}^3\times(0,t)$, and then using the continuity equation \eqref{ma1}$_1$ and integration by parts, we obtain that
\begin{align}
&\frac{1}{2}\int_{\mathbb{R}^3}
\big(|\sqrt{\rho^\epsilon}u^\epsilon|^2+|\eta^\epsilon-1|^2
+|\tau^\epsilon|^2+2\Pi^{\epsilon}\big)\;\mathrm{d}x\notag\\
&+\int_0^t\int_{\mathbb{R}^3}\big(\mu_1|\nabla u^\epsilon|^2+\mu_2|\mathrm{div} u^\epsilon|^2+\nu|\nabla \eta^\epsilon|^2+\frac{A_0}{2}|\tau^\epsilon|^2+\nu|\nabla \tau^\epsilon|^2\big) \;\mathrm{d}x\mathrm{d}s\label{s43}\\
=&\frac{1}{2}\int_{\mathbb{R}^3}
\big(|\sqrt{\rho^\epsilon_0}u^\epsilon_0|^2+|\eta^\epsilon_0-1|^2
+|\tau^\epsilon_0|^2+2\Pi^{\epsilon}_0\big)\;\mathrm{d}x-\beta(L-1)\int_0^t\int_{\mathbb{R}^3}u^\epsilon\cdot \nabla\eta^\epsilon\;\mathrm{d}x\notag\\
&-2\bar{\mathfrak{z}}\int_0^t\int_{\mathbb{R}^3}\eta^\epsilon u^\epsilon\cdot\nabla\eta^\epsilon \;\mathrm{d}x\mathrm{d}s
+\frac{\beta}{k}\int_0^t\int_{\mathbb{R}^3} \mathrm{div}\tau^\epsilon \cdot u^\epsilon\;\mathrm{d}x\mathrm{d}s-\int_0^t\int_{\mathbb{R}^3} \mathrm{div}u^\epsilon(\eta^\epsilon-1)\,\mathrm{d}x\mathrm{d}s\notag\\
&+\int_0^t\int_{\mathbb{R}^3}(\eta^\epsilon-1)u^\epsilon\cdot\nabla(\eta^\epsilon-1)\,\mathrm{d}x\mathrm{d}s
+\int_0^t\int_{\mathbb{R}^3} (u^\epsilon\cdot\nabla)\tau^\epsilon:\tau^\epsilon\;\mathrm{d}x\mathrm{d}s\notag\\
&+\int_0^t\int_{\mathbb{R}^3} (\nabla u^\epsilon\tau^\epsilon+\tau^\epsilon\nabla^T u^\epsilon):\tau^\epsilon\;\mathrm{d}x\mathrm{d}s
+k\int_0^t\int_{\mathbb{R}^3}\eta^\epsilon (\nabla u^\epsilon+\nabla^T u^\epsilon):\tau^\epsilon\;\mathrm{d}x\mathrm{d}s\notag.
\end{align}
Meanwhile, multiplying the equations \eqref{limit1} by $(u, \eta-1, \tau)$, integrating the resulting equations over $\mathbb{R}^3\times(0,t)$, and using integration by parts, we have
\begin{align}
&\frac{1}{2}\int_{\mathbb{R}^3}
\big(|u|^2+|\eta-1|^2
+|\tau|^2\big)\;\mathrm{d}x
+\int_0^t\int_{\mathbb{R}^3}\big(\mu_1|\nabla u|^2+\nu|\nabla \eta|^2+\frac{A_0}{2}|\tau|^2+\nu|\nabla \tau|^2\big) \;\mathrm{d}x\mathrm{d}s\label{s44}\\
=&\frac{1}{2}\int_{\mathbb{R}^3}
\big(|u_0|^2+|\eta_0-1|^2
+|\tau_0|^2\big)\;\mathrm{d}x+\frac{\beta}{k}\int_0^t\int_{\mathbb{R}^3} \mathrm{div}\tau \cdot u\;\mathrm{d}x\mathrm{d}s\notag\\
&+\int_0^t\int_{\mathbb{R}^3}(\eta-1)u\cdot\nabla(\eta-1)\,\mathrm{d}x\mathrm{d}s-\int_0^t\int_{\mathbb{R}^3} \mathrm{div}u(\eta-1)\,\mathrm{d}x\mathrm{d}s+\int_0^t\int_{\mathbb{R}^3} (u\cdot\nabla)\tau:\tau\;\mathrm{d}x\mathrm{d}s\notag\\
&+\int_0^t\int_{\mathbb{R}^3} (\nabla u\tau+\tau\nabla^T u):\tau\;\mathrm{d}x\mathrm{d}s
+k\int_0^t\int_{\mathbb{R}^3}\eta (\nabla u+\nabla^T u):\tau\;\mathrm{d}x\mathrm{d}s\notag.
\end{align}
On the other hand, by the basic principles of calculus, combining equations \eqref{ma1}$_2$ with \eqref{limit1}$_1$ and using integration by parts, we can get the following equality
\begin{align}
&\int_{\mathbb{R}^3}
\rho^\epsilon u^\epsilon\cdot u\;\mathrm{d}x
-\int_{\mathbb{R}^3}
\rho^\epsilon_0 u^\epsilon_0\cdot u_0\;\mathrm{d}x\notag\\
=&\int_0^t\int_{\mathbb{R}^3}
(\rho^\epsilon u^\epsilon)_t\cdot u\;\mathrm{d}x\mathrm{d}s
+\int_0^t\int_{\mathbb{R}^3}
\rho^\epsilon u^\epsilon\cdot u_t\;\mathrm{d}x\mathrm{d}s\label{s45}\\
=&\int_0^t\int_{\mathbb{R}^3}
\rho^\epsilon u^\epsilon\otimes u^\epsilon:\nabla u
\;\mathrm{d}x\mathrm{d}s
+\mu_1\int_0^t\int_{\mathbb{R}^3}\Delta u^\epsilon\cdot u\;\mathrm{d}x\mathrm{d}s
+\frac{\beta}{k}\int_0^t\int_{\mathbb{R}^3}\mathrm{div}\tau^\epsilon\cdot u\;\mathrm{d}x\mathrm{d}s\notag\\
&-\int_0^t\int_{\mathbb{R}^3}\rho^\epsilon u^\epsilon\otimes u:\nabla u\;\mathrm{d}x\mathrm{d}s
-\int_0^t\int_{\mathbb{R}^3}\rho^\epsilon u^\epsilon\cdot\nabla\pi\;\mathrm{d}x\mathrm{d}s
+\mu_1\int_0^t\int_{\mathbb{R}^3}\rho^\epsilon u^\epsilon\cdot\Delta u\;\mathrm{d}x\mathrm{d}s\notag\\
&+\frac{\beta}{k}\int_0^t\int_{\mathbb{R}^3}\rho^\epsilon u^\epsilon\cdot\mathrm{div}\tau\;\mathrm{d}x\mathrm{d}s\notag.
\end{align}
Further, \eqref{s45} can be rewritten as
\begin{align}
&-\int_{\mathbb{R}^3}
\sqrt{\rho^\epsilon} u^\epsilon\cdot u\;\mathrm{d}x
-2\mu_1\int_0^t\int_{\mathbb{R}^3}\nabla u^\epsilon:\nabla u\;\mathrm{d}x\mathrm{d}s\notag\\
=&-\int_{\mathbb{R}^3}
\rho^\epsilon_0 u^\epsilon_0\cdot u_0\;\mathrm{d}x
+\int_{\mathbb{R}^3}(\rho^\epsilon-\sqrt{\rho^\epsilon})u^\epsilon\cdot u\;\mathrm{d}x+\int_0^t\int_{\mathbb{R}^3}
\rho^\epsilon u^\epsilon\otimes (u-u^\epsilon):\nabla u
\;\mathrm{d}x\mathrm{d}s\label{s46}\\
&-\mu_1\int_0^t\int_{\mathbb{R}^3}(\rho^\epsilon u^\epsilon-u^\epsilon)\cdot\Delta u\;\mathrm{d}x\mathrm{d}s
+\int_0^t\int_{\mathbb{R}^3}\rho^\epsilon u^\epsilon\cdot\nabla\pi\;\mathrm{d}x\mathrm{d}s\notag\\
&-\frac{\beta}{k}\int_0^t\int_{\mathbb{R}^3}\mathrm{div}\tau\cdot u^\epsilon\;\mathrm{d}x\mathrm{d}s-\frac{\beta}{k}\int_0^t\int_{\mathbb{R}^3}\mathrm{div}\tau^\epsilon\cdot u\;\mathrm{d}x\mathrm{d}s
-\frac{\beta}{k}\int_0^t\int_{\mathbb{R}^3}\mathrm{div}\tau\cdot (\rho^\epsilon-1)u^\epsilon\;\mathrm{d}x\mathrm{d}s\notag.
\end{align}
Similarly, based on the basic principles of calculus and integration by parts, we deduce the following inequality from equations \eqref{ma1}$_3$ and \eqref{limit1}$_2$:
\begin{align}
&-\int_{\mathbb{R}^3}
(\eta^\epsilon-1)\cdot (\eta-1)\,\mathrm{d}x
+\int_{\mathbb{R}^3}
(\eta^\epsilon_0-1) \cdot (\eta_0-1)\,\mathrm{d}x\notag\\
=&-\int_0^t\int_{\mathbb{R}^3}
(\eta^\epsilon-1)_t\cdot (\eta-1)\,\mathrm{d}x\mathrm{d}s
-\int_0^t\int_{\mathbb{R}^3}
(\eta^\epsilon-1) \cdot (\eta-1)_t\,\mathrm{d}x\mathrm{d}s\label{s47}\\
=&\nu\int_0^t\int_{\mathbb{R}^3}
\nabla\eta^\epsilon\cdot\nabla \eta\,\mathrm{d}x\mathrm{d}s
-\int_0^t\int_{\mathbb{R}^3}
(\eta^\epsilon-1) u^\epsilon\cdot\nabla( \eta-1)\,\mathrm{d}x\mathrm{d}s\notag\\
&+\nu\int_0^t\int_{\mathbb{R}^3}
\nabla\eta\cdot\nabla \eta^\epsilon\,\mathrm{d}x\mathrm{d}s
-\int_0^t\int_{\mathbb{R}^3}
(\eta-1) u\cdot\nabla (\eta^\epsilon-1)\,\mathrm{d}x\mathrm{d}s\notag\\
&+\int_0^t\int_{\mathbb{R}^3}
\mathrm{div}u^\epsilon(\eta-1)\,\mathrm{d}x\mathrm{d}s+\int_0^t\int_{\mathbb{R}^3}
\mathrm{div}u(\eta^\epsilon-1)\,\mathrm{d}x\mathrm{d}s\notag,
\end{align}
as well as the following inequality from equations \eqref{ma1}$_4$ and\eqref{limit1}$_3$:
\begin{align}
&-\int_{\mathbb{R}^3}
\tau^\epsilon: \tau\;\mathrm{d}x
+\int_{\mathbb{R}^3}
\tau^\epsilon_0: \tau_0\;\mathrm{d}x\notag\\
=&-\int_0^t\int_{\mathbb{R}^3}
\tau^\epsilon_t: \tau\;\mathrm{d}x\mathrm{d}s
-\int_0^t\int_{\mathbb{R}^3}
\tau^\epsilon : \tau_t\;\mathrm{d}x\mathrm{d}s\label{s48}\\
=&\frac{A_0}{2}\int_0^t\int_{\mathbb{R}^3}
\tau^\epsilon:\tau\;\mathrm{d}x\mathrm{d}s
+\nu\int_0^t\int_{\mathbb{R}^3}
\nabla\tau^\epsilon:\nabla \tau\;\mathrm{d}x\mathrm{d}s
-\int_0^t\int_{\mathbb{R}^3}
\tau^\epsilon u^\epsilon\cdot\nabla \tau\;\mathrm{d}x\mathrm{d}s\notag\\
&-\int_0^t\int_{\mathbb{R}^3}
(\nabla u^\epsilon\tau^\epsilon+\tau^\epsilon\nabla^T u^\epsilon): \tau\;\mathrm{d}x\mathrm{d}s
-k\int_0^t\int_{\mathbb{R}^3}
\eta^\epsilon(\nabla u^\epsilon+\nabla^T u^\epsilon): \tau\;\mathrm{d}x\mathrm{d}s\notag\\
&+\frac{A_0}{2}\int_0^t\int_{\mathbb{R}^3}
\tau:\tau^\epsilon\;\mathrm{d}x\mathrm{d}s
+\nu\int_0^t\int_{\mathbb{R}^3}
\nabla\tau:\nabla \tau^\epsilon\;\mathrm{d}x\mathrm{d}s
-\int_0^t\int_{\mathbb{R}^3}
\tau u\cdot\nabla \tau^\epsilon\;\mathrm{d}x\mathrm{d}s\notag\\
&-\int_0^t\int_{\mathbb{R}^3}
(\nabla u\tau+\tau\nabla^T u): \tau^\epsilon\;\mathrm{d}x\mathrm{d}s
-k\int_0^t\int_{\mathbb{R}^3}
\eta(\nabla u+\nabla^T u): \tau^\epsilon\;\mathrm{d}x\mathrm{d}s\notag.
\end{align}
A reasonable combination of \eqref{s43}, \eqref{s44} and \eqref{s46}-\eqref{s48} yields that
\begin{align}
&\frac{1}{2}\int_{\mathbb{R}^3}
\big(|\sqrt{\rho^\epsilon}u^\epsilon-u|^2+|\eta^\epsilon-\eta|^2
+|\tau^\epsilon-\tau|^2+2\Pi^{\epsilon}\big)\;\mathrm{d}x\notag\\
&+\int_0^t\int_{\mathbb{R}^3}\big(\mu_1|\nabla u^\epsilon-\nabla u|^2+\mu_2|\mathrm{div} u^\epsilon|^2+\nu|\nabla \eta^\epsilon-\nabla \eta|^2+\frac{A_0}{2}|\tau^\epsilon-\tau|^2+\nu|\nabla \tau^\epsilon-\nabla \tau|^2\big) \;\mathrm{d}x\mathrm{d}s\label{s49}\\
=&\frac{1}{2}\int_{\mathbb{R}^3}
\left(|\sqrt{\rho^\epsilon_0}u^\epsilon_0-u_0|^2+|\eta^\epsilon_0-\eta_0|^2
+|\tau^\epsilon_0-\tau_0|^2+2\Pi^{\epsilon}_0\right)\,\mathrm{d}x
-\int_{\mathbb{R}^3}\sqrt{\rho^\epsilon_0}(\sqrt{\rho^\epsilon_0}-1)u^\epsilon_0\cdot u_0\,\mathrm{d}x\notag\\
&+\int_{\mathbb{R}^3}(\rho^\epsilon-\sqrt{\rho^\epsilon})u^\epsilon\cdot u\,\mathrm{d}x-\mu_1\int_0^t\int_{\mathbb{R}^3}(\rho^\epsilon -1)u^\epsilon\cdot\Delta u\,\mathrm{d}x\mathrm{d}s
-\beta(L-1)\int_0^t\int_{\mathbb{R}^3}u^\epsilon\cdot \nabla\eta^\epsilon\,\mathrm{d}x\mathrm{d}s\notag\\
&-2\bar{\mathfrak{z}}\int_0^t\int_{\mathbb{R}^3}\eta^\epsilon u^\epsilon\cdot\nabla\eta^\epsilon \,\mathrm{d}x\mathrm{d}s+\int_0^t\int_{\mathbb{R}^3}\rho^\epsilon u^\epsilon\cdot\nabla\pi\,\mathrm{d}x\mathrm{d}s
-\frac{\beta}{k}\int_0^t\int_{\mathbb{R}^3}\mathrm{div}\tau\cdot (\rho^\epsilon-1)u^\epsilon\,\mathrm{d}x\mathrm{d}s\notag\\
&+\frac{\beta}{k}\int_0^t\int_{\mathbb{R}^3}(\mathrm{div}\tau^\epsilon-\mathrm{div}\tau)\cdot (u^\epsilon-u)\,\mathrm{d}x\mathrm{d}s+\int_0^t\int_{\mathbb{R}^3}(\eta-\eta^\epsilon)(\mathrm{div}u^\epsilon-\mathrm{div}u)\,\mathrm{d}x\mathrm{d}s\notag\\
&+\int_0^t\int_{\mathbb{R}^3}[(\eta^\epsilon-1)u^\epsilon-(\eta-1) u]\cdot(\nabla\eta^\epsilon-\nabla\eta)\,\mathrm{d}x\mathrm{d}s
\notag\\
&
+k\int_0^t\int_{\mathbb{R}^3}[\eta^\epsilon (\nabla u^\epsilon+\nabla^T u^\epsilon)-\eta (\nabla u+\nabla^T u)]:(\tau^\epsilon-\tau)\,\mathrm{d}x\mathrm{d}s\notag\\[1mm]
&+\int_0^t\int_{\mathbb{R}^3}
\rho^\epsilon u^\epsilon\otimes (u-u^\epsilon):\nabla u
\,\mathrm{d}x\mathrm{d}s+\int_0^t\int_{\mathbb{R}^3} (\tau^\epsilon u^\epsilon-\tau u)\cdot\nabla(\tau^\epsilon-\tau)\,\mathrm{d}x\mathrm{d}s\notag\\[1mm]
&+\int_0^t\int_{\mathbb{R}^3}
[(\nabla u^\epsilon\tau^\epsilon+\tau^\epsilon\nabla^T u^\epsilon)-(\nabla u\tau+\tau\nabla^T u)]: (\tau^\epsilon-\tau)\,\mathrm{d}x\mathrm{d}s
:=\sum\limits_{i=1}^{15}II_i.\notag
\end{align}
In what follows, we estimate the terms $II_1$-$II_{15}$. Observing that the terms $II_2$, $II_3$, $II_4$ and $II_{13}$ are the same as the terms stated in Section 4 of \cite{GJLLT}, we can directly get
\begin{align}
II_2+II_3+II_4+II_{13}
\leq& C\epsilon^{1+\frac{\alpha_0}{2}}+C_{\theta_4}(1+T)\epsilon^2
+\theta_4\langle\Pi^\epsilon,1\rangle\label{s50}\\
&+C\int_0^t\langle\Pi^\epsilon,1\rangle\;\mathrm{d}s
+C\int_0^t\|\sqrt{\rho^\epsilon}  u ^\epsilon- u \|_{L^2}^2\;\mathrm{d}s\notag,
\end{align}
where $\theta_4>0$ is a small constant independent of$\epsilon$. For $II_1$ and $II_8$, we can use \eqref{clc-02}, \eqref{clc-09}, \eqref{clc-10}, Lemma \ref{lemma1}, H\"older inequality, Sobolev inequality and Cauchy inequality to obtain
\begin{align}
II_1+II_8
\lesssim& \epsilon^{\alpha_0}+
\int_0^t\int_{\mathbb{R}^3}|\nabla\tau|\, |\rho^\epsilon-1|\,|u^\epsilon|\;\mathrm{d}x\mathrm{d}s\notag\\
\lesssim&
\epsilon^{\alpha_0}+\int_0^t\|\sqrt{\rho^\epsilon}-1\|_{L^2}\|\sqrt{\rho^\epsilon}+1\|_{L^\infty}
\|u^\epsilon\|_{L^\infty}\|\nabla\tau\|_{L^2}\;\mathrm{d}s\label{s51}\\
\lesssim& \epsilon^{\alpha_0}+ \int_0^t\epsilon\langle\Pi^\epsilon,1\rangle^\frac{1}{2}\;\mathrm{d}s\notag\\
\lesssim& \epsilon^{\alpha_0}+ T\epsilon^2
+\int_0^t\langle\Pi^\epsilon,1\rangle\;\mathrm{d}s\notag.
\end{align}
Recalling that $\rho^\epsilon_t=-\mathrm{div}(\rho^\epsilon u^\epsilon)$ and $\pi=\pi_1+\beta(L-1)\eta+\bar{\mathfrak{z}}(\eta)^{2}$, we use integration by parts to rewrite $II_i\ (i=5,6,7)$ as
\begin{align}
\sum_{i=5}^{7}II_i=&\beta(L-1)\int_0^t\int_{\mathbb{R}^3}(\rho^\epsilon u^\epsilon\cdot\nabla\eta-u^\epsilon\cdot \nabla\eta^\epsilon)\;\mathrm{d}x\mathrm{d}s
+\bar{\mathfrak{z}}\int_0^t\int_{\mathbb{R}^3}[\rho^\epsilon u^\epsilon\cdot\nabla(\eta^{2})- u^\epsilon\cdot\nabla(\eta^{\epsilon2})] \;\mathrm{d}x\mathrm{d}s\notag\\
&+\int_0^t\int_{\mathbb{R}^3}\rho^\epsilon u^\epsilon\cdot\nabla\pi_1\;\mathrm{d}x\mathrm{d}s\notag\\
=&\beta(L-1)\int_0^t\int_{\mathbb{R}^3}[(\rho^\epsilon-1) u^\epsilon\cdot \nabla\eta^\epsilon+ \rho^\epsilon u^\epsilon\cdot\nabla(\eta-\eta^\epsilon)]\;\mathrm{d}x\mathrm{d}s\notag\\
&+\bar{\mathfrak{z}}\int_0^t\int_{\mathbb{R}^3}[(\rho^\epsilon -1)u^\epsilon\cdot\nabla(\eta^{\epsilon2})+ \rho^\epsilon u^\epsilon\cdot\nabla(\eta^{2}-\eta^{\epsilon2})] \;\mathrm{d}x\mathrm{d}s
+\int_0^t\int_{\mathbb{R}^3}\rho^\epsilon u^\epsilon\cdot\nabla\pi_1\;\mathrm{d}x\mathrm{d}s\notag\\
=&\int_0^t\int_{\mathbb{R}^3}(\sqrt{\rho^\epsilon}-1)(\sqrt{\rho^\epsilon}+1) [\beta(L-1)u^\epsilon\cdot \nabla\eta^\epsilon+\bar{\mathfrak{z}}u^\epsilon\cdot\nabla(\eta^{\epsilon2})]\;\mathrm{d}x\mathrm{d}s\notag\\
&+\int_0^t\int_{\mathbb{R}^3}\rho^\epsilon_t [\beta(L-1)(\eta-\eta^\epsilon)+\bar{\mathfrak{z}} (\eta^{2}-\eta^{\epsilon2})+\pi_1] \;\mathrm{d}x\mathrm{d}s\notag.
\end{align}
Further, similar to \eqref{s51}, we can deduce
\begin{align}
&\int_0^t\int_{\mathbb{R}^3}\rho^\epsilon_t [\beta(L-1)(\eta-\eta^\epsilon)+\bar{\mathfrak{z}} (\eta^{2}-\eta^{\epsilon2})+\pi_1] \;\mathrm{d}x\mathrm{d}s\notag\\
=&\int_0^t\int_{\mathbb{R}^3}(\rho^\epsilon-1)_t [\beta(L-1)(\eta-\eta^\epsilon)+\bar{\mathfrak{z}} (\eta^{2}-\eta^{\epsilon2})+\pi_1] \;\mathrm{d}x\mathrm{d}s\label{s52}\\
=&\int_{\mathbb{R}^3}(\sqrt{\rho^\epsilon}-1)(\sqrt{\rho^\epsilon}+1) [\beta(L-1)(\eta-\eta^\epsilon)+\bar{\mathfrak{z}} (\eta^{2}-\eta^{\epsilon2})+\pi_1] \;\mathrm{d}x\notag\\
&-\int_{\mathbb{R}^3}(\sqrt{\rho^\epsilon_0}-1)(\sqrt{\rho^\epsilon_0}+1) [\beta(L-1)(\eta_0-\eta^\epsilon_0)+\bar{\mathfrak{z}} (\eta^{2}_0-\eta^{\epsilon2}_0)+\pi_1(0)] \;\mathrm{d}x\notag\\
&-\int_0^t\int_{\mathbb{R}^3}(\sqrt{\rho^\epsilon}-1)(\sqrt{\rho^\epsilon}+1) [\beta(L-1)(\eta-\eta^\epsilon)_t+\bar{\mathfrak{z}} (\eta^{2}-\eta^{\epsilon2})_t+(\pi_1)_t] \;\mathrm{d}x\mathrm{d}s\notag\\
\leq&C_{\theta_4}(1+T)\epsilon^2+C\epsilon^{1+\frac{\alpha_0}{2}}
+\theta_4\langle\Pi^\epsilon,1\rangle
+C\int_0^t\langle\Pi^\epsilon,1\rangle\;\mathrm{d}s\notag,
\end{align}
and
\begin{align}
\int_0^t\int_{\mathbb{R}^3}(\sqrt{\rho^\epsilon}-1)(\sqrt{\rho^\epsilon}+1) [\beta(L-1)u^\epsilon\cdot \nabla\eta^\epsilon+\bar{\mathfrak{z}}u^\epsilon\cdot\nabla(\eta^{\epsilon2})]\;\mathrm{d}x\mathrm{d}s
\lesssim T\epsilon^2
+\int_0^t\langle\Pi^\epsilon,1\rangle\;\mathrm{d}s.\label{s53}
\end{align}
Combining \eqref{s52} with \eqref{s53}, we find that
\begin{align}
II_5+II_6+II_7\leq C_{\theta_4}(1+T)\epsilon^2+C\epsilon^{1+\frac{\alpha_0}{2}}
+\theta_4\langle\Pi^\epsilon,1\rangle
+C\int_0^t\langle\Pi^\epsilon,1\rangle\;\mathrm{d}s.\label{s54}
\end{align}
For $II_9$, $II_{10}$, $II_{11}$ and $II_{14}$, by the virtue of \eqref{clc-02}, \eqref{clc-09}, H\"older inequality, Sobolev inequality, Cauchy inequality and Lemma \ref{lemma1}, we have
\begin{align}
II_9=&\frac{\beta}{k}\int_0^t\int_{\mathbb{R}^3}
(\mathrm{div}\tau^\epsilon-\mathrm{div}\tau)\cdot u^\epsilon(1-\sqrt{\rho^\epsilon})\;\mathrm{d}x\mathrm{d}s
+\frac{\beta}{k}\int_0^t\int_{\mathbb{R}^3}
(\mathrm{div}\tau^\epsilon-\mathrm{div}\tau)\cdot (\sqrt{\rho^\epsilon}u^\epsilon-u)\;\mathrm{d}x\mathrm{d}s\notag\\
\leq&\theta_4\int_0^t\|\mathrm{div}\tau^\epsilon-\mathrm{div}\tau\|_{L^2}^2\;\mathrm{d}s
+C_{\theta_4}\int_0^t\|\sqrt{\rho^\epsilon}-1\|_{L^2}^2\;\mathrm{d}s
+C_{\theta_4}\int_0^t\|\sqrt{\rho^\epsilon}u^\epsilon-u\|_{L^2}^2\;\mathrm{d}s\label{s55}\\
\leq&\theta_4\int_0^t\|\nabla\tau^\epsilon-\nabla\tau\|_{L^2}^2\;\mathrm{d}s
+C_{\theta_4}\int_0^t\langle\Pi^\epsilon,1\rangle\;\mathrm{d}s
+C_{\theta_4}\int_0^t\|\sqrt{\rho^\epsilon}u^\epsilon-u\|_{L^2}^2\;\mathrm{d}s\notag,
\end{align}
and
\begin{align}
&II_{10}+II_{11}+II_{14}\notag\\
=&\int_0^t\int_{\mathbb{R}^3}(\eta-\eta^\epsilon)(\mathrm{div}u^\epsilon-\mathrm{div}u)\,\mathrm{d}x\mathrm{d}s+\int_0^t\int_{\mathbb{R}^3}(\eta^\epsilon -\eta)u\cdot(\nabla\eta^\epsilon-\nabla\eta)\,\mathrm{d}x\mathrm{d}s\notag\\
&+\int_0^t\int_{\mathbb{R}^3}(\eta^\epsilon -1)[(1-\sqrt{\rho^\epsilon})u^\epsilon
+(\sqrt{\rho^\epsilon}u^\epsilon-u)]\cdot(\nabla\eta^\epsilon-\nabla\eta)\,\mathrm{d}x\mathrm{d}s\label{s56}\\
&+\int_0^t\int_{\mathbb{R}^3} [(\tau^\epsilon -\tau)u^\epsilon+\tau(1-\sqrt{\rho^\epsilon})u^\epsilon
+\tau(\sqrt{\rho^\epsilon}u^\epsilon-u)]\cdot(\nabla\tau^\epsilon-\nabla\tau)\,\mathrm{d}x\mathrm{d}s
\notag\\
\leq&\theta_4\Big(\int_0^t\|\nabla u^\epsilon-\nabla u\|_{L^2}^2\;\mathrm{d}s+\int_0^t\|\nabla\eta^\epsilon-\nabla\eta\|_{L^2}^2\;\mathrm{d}s+
\int_0^t\|\nabla\tau^\epsilon-\nabla\tau\|_{L^2}^2\;\mathrm{d}s\Big)+C_{\theta_4}\epsilon^2 T
\notag\\
&+C_{\theta_4}\Big(\int_0^t\|\eta^\epsilon-\eta\|_{L^2}^2\;\mathrm{d}s
+\int_0^t\|\tau^\epsilon-\tau\|_{L^2}^2\;\mathrm{d}s
+\int_0^t\|\sqrt{\rho^\epsilon}u^\epsilon-u\|_{L^2}^2\;\mathrm{d}s
+\int_0^t\langle\Pi^\epsilon,1\rangle\;\mathrm{d}s\Big).\notag
\end{align}
Similarly, $II_{12}$ and $II_{15}$ can be controlled as follows
\begin{align}\label{s57}
II_{12}+II_{15}\leq \theta_4\int_0^t\|\nabla u^\epsilon-\nabla u\|_{L^2}^2\;\mathrm{d}s
+C_{\theta_4}\int_0^t\|\tau^\epsilon-\tau\|_{L^2}^2\;\mathrm{d}s
+C_{\theta_4}\int_0^t\|\eta^\epsilon-\eta\|_{L^2}^2\;\mathrm{d}s.
\end{align}
In conclusion, substituting \eqref{s50}, \eqref{s51}, \eqref{s54}-\eqref{s57} into \eqref{s49}, choosing $\theta_4$ sufficiently small, we obtain that
\begin{align}
&\int_{\mathbb{R}^3}
\big(|\sqrt{\rho^\epsilon}u^\epsilon-u|^2+|\eta^\epsilon-\eta|^2
+|\tau^\epsilon-\tau|^2+2\Pi^{\epsilon}\big)\;\mathrm{d}x\notag\\
&+\int_0^t\int_{\mathbb{R}^3}\big(\mu_1|\nabla u^\epsilon-\nabla u|^2+\mu_2|\mathrm{div} u^\epsilon|^2+\nu|\nabla \eta^\epsilon-\nabla \eta|^2+\frac{A_0}{2}|\tau^\epsilon-\tau|^2+\nu|\nabla \tau^\epsilon-\nabla \tau|^2\big) \;\mathrm{d}x\mathrm{d}s\label{s63}\\
\lesssim &(1+T)\epsilon^{\beta_0}
+\int_0^t\big(\|\sqrt{\rho^\epsilon} u^\epsilon-u\|_{L^2}^2+\|\eta^\epsilon-\eta\|_{L^2}^2
+\|\tau^\epsilon-\tau\|_{L^2}^2
+\langle\Pi^\epsilon,1\rangle\big)\;\mathrm{d}s\notag,
\end{align}
where $\beta_0$ is defined in Theorem \ref{Theorem2}. Based on \eqref{s63} and Gronwall's inequality, we can easily get \eqref{clc-11} and finish the proof of Theorem \ref{Theorem2}.
\endpf

\section*{Acknowledgement}
The authors would like to thank Prof. Huanyao Wen for valuable discussions and suggestions. This paper was supported by the National Natural Science Foundation of China $\#$12071152 and the Natural Science Foundation of Guangdong Province $\#$ 2022A1515012112.

\end{document}